\documentclass[a4paper,10pt]{article}

\usepackage{hyperref}
\usepackage[utf8]{inputenc}
\usepackage{amssymb}
\usepackage{amsmath}
\usepackage{graphicx}
\usepackage{xcolor}
\usepackage{natbib}

\usepackage{algorithm}
\usepackage{algpseudocode}
\usepackage{algorithmicx}

\usepackage[left=2cm,top=2.5cm,right=2cm,bottom=2.5cm]{geometry}

\hypersetup{colorlinks,linkcolor= {blue},citecolor= {blue}}

\graphicspath{{Figures/}}

\title{3D fictitious wave domain CSEM inversion by adjoint source estimation}

\author{Pengliang Yang$^1$\\
  $^1$ School of Mathematics, Harbin Institute of Technology, 150001, Harbin, China\\
  E-mail: ypl.2100@gmail.com
}

\begin{document}

\maketitle

\begin{abstract}
  Marine controlled-source electromagnetic (CSEM) method has proved its potential in detecting
  highly resistive hydrocarbon bearing formations.  A novel frequency domain CSEM inversion approach using fictitious wave domain time stepping modelling is presented.
  Using Lagrangian-based adjoint state method, the inversion gradient with respect to resistivity  can be computed by the product between the forward and adjoint fields. Simulation of the adjoint field using the same modelling engine is challenging as it requires time domain adjoint source time functions while only a few discrete frequencies of the data residual are available for the inversion. A regularized linear inverse problem is formulated in order to estimate a long time series from very few frequency samples. It can then be solved using linear optimization technique, yielding a matrix-free implementation. Instead of computing adjoint source time function one by one at each receiver location, a basis function implementation has been developed such that the inverse problem can be solved only once and reused every time to construct all time-domain adjoint sources. The method allows computing all frequencies of the EM fields  in one go without heavy memory and computational overhead, making efficient 3D CSEM inversion feasible. Numerical examples are employed to demonstrate the application of our method.
\end{abstract}

\section{Introduction}

Controlled-source electromagnetic (CSEM) method is a well established technology to do geophysical exploration \citep{chave1982on,constable1986offshore}. It can also be configured for air-borne \citep{yin2015review} and cross-well \citep{alumbaugh1997three} geometries.  Land CSEM has been a commonplace to find mineral deposit \citep{ward1988electromagnetic,zhdanov1994geoelectrical,grayver2014geo}. As a complement to seismic measurement, marine CSEM have been successfully applied to detect hydrocarbon bearing formations in oil and gas industry 
\citep{eidesmo2002sbl,ellingsrud2002remote,constable2007introduction,macgregor2014mcsem}, using vessel-towed dipole source over an array of receivers deployed on the seabed (hence the name seabed logging). The key to these applications is the good sensitivity of the electromagnetic signals to distinguish high contrast in resistivity between saline-filled rocks and hydrocarbons \citep{abubakar20082}. This property makes CSEM an ideal tool for de-risking \citep{macgregor2007derisking} in reservoir exploration prior to drilling, as well as 4D monitoring during the production \citep{shantsev2002timelapse}. A topical review on marine CSEM inversion has been given by \citet{constable2010}.

CSEM inversion, also known as resistivity tomography, has been a standard technique to translate CSEM data into electrical properties of the subsurface. To determine the resistivity model, a nonlinear inverse problem is then formulated to iteratively minimize the difference between the observed EM data and the synthetic data derived from numerical modelling, using different optimization schemes, such as nonlinear conjugate gradient method \citep{gribenko2007rigorous,commer2008new}, quasi-Newton l-BFGS \citep{Plessix_2008_RIC,schwarzbach2013gji}, and Gauss-Newton method \citep{constable1987occam,abubakar20082,zaslavsky2013large}.

The most computation intensive part of CSEM inversion is the numerical simulation of 3D electromagnetic field. One can consider the time-domain finite-difference method \citep{oristaglio1984diffusion,wang1993finite,Taflove_2005_CEF}, the frequency-domain finite-difference method \citep{newman1995frequency,smith1996conservative1,mulder2006multigrid,streich20093d}, and the frequency-domain finite-element method \citep{li20072d,da2012finite,key2016mare2dem,rochlitz2019custem}.

A plethora of CSEM modelling tools directly solve the linear system based on the frequency domain solution of the discretized Maxwell equation.
This avoids the high computational cost by direct discretization of time-domain diffusive Maxwell equation involving extremely large number of time steps dictated by the restrictive stability condition. Both the frequency domain finite difference method and finite element method formulate Maxwell equation as a matrix-based linear system, which may be solved using direct \citep{streich20093d} or iterative \citep{smith1996conservative2,mulder2006multigrid,puzyrev2013parallel} solvers.
Modelling by direct solver is very attractive for multi-source CSEM problems \citep{streich20093d}, as the matrix system can be factorized only once and reused for all sources by forward-backward substitution. Since the inversion of the large sparse matrix for 3D problems involves huge amount of memory resources, direct methods may easily go beyond the memory capacity of a desktop computer. Iterative solvers require much less memory storage but may be time-consuming and difficult to converge due to ill-conditioning of the discretized Helmholtz matrix when the conductivity/resistivity model is highly heterogeneous or the modelling grid is severely stretched \citep{mulder2006multigrid}.
While the frequency domain finite element method is flexible to address complex model geometries thanks to the pre-computed meshes, the meshing in 3D geometries using tetrahedral and hexahedral mesh itself is a challenging and time consuming task.

In the time-domain, the fields are updated at each time step, allowing the modelling over the same memory units. The frequency domain fields can be integrated on the fly during time-stepping, such that multiple frequencies can be extracted from the same simulation. This avoids repeating several times of the modelling for different frequencies in frequency domain methods.
It motivates \citet{Maao_2007_FFT} to propose a modified wave domain approach to significantly speed up the computation of diffusive electromagnetic modelling. The method has been adapted in \citet{Storen_2008_Gradient} to efficiently compute the inversion gradient for 3D industrial scale applications. The modified wave domain approach by \citet{Maao_2007_FFT} has an attenuation/diffusive term. Inspired by the fictitious wave domain approach initially proposed in \citet{lee1989new},  \citet{Mittet_2010_HFD} transformed the diffusion domain Maxwell equation into a pure wave domain. 

A straightforward time-domain discretization of the diffusive Maxwell equation leads to the temporal sampling proportional to the square of grid spacing, due to the requirement of the stability condition ($\Delta t \propto \Delta x^2$). This means taking half the spatial sampling will lead to a quadratic increase of the number of time steps for a simulation of the same duration. Transforming the diffusive domain into wave domain using fictitious wave domain approach allows a linear increase of time step with the use of finer grid spacing ($\Delta t\propto \Delta x$). This yields a highly efficient scheme to compute the same frequency domain EM field, as it significantly reduces the required number of time steps. These advantages inspire us to adopt the fictitious wave domain method for CSEM inversion.

The major contribution of this paper is to develop an efficient frequency domain CSEM inversion scheme based on the fictitious wave domain modelling. The gradient of the misfit functional for frequency domain CSEM inversion requires both a forward and an adjoint field. Efficient computation of the frequency domain adjoint field by time domain modelling manifests itself as a major challenge in this development as it requires the time domain adjoint source time functions at all receiver locations which are unavailable due to the formulation of CSEM inversion in frequency domain. The key novelty of our approach is to estimate the time-domain adjoint source based on only a few discrete frequency samples of the data residual. To do so, we formulate a regularized linear inverse problem, and then present an efficient, matrix-free implementation using very few basis functions. The linear inverse problem can be solved only once and the resulting solution can be reused to derive the adjoint source time functions at all receiver locations. This novel development delivers an efficient 3D CSEM inversion methodology in fictitious wave domain. We finally apply our method to two examples for numerical demonstration.

\section{CSEM modelling in fictitious wave domain}\label{sec:forward}

Within a three dimensional space $X$, the diffusive Maxwell equations  are governed  in the frequency domain
\begin{equation}\label{eq:freq}
  \begin{cases}
    \nabla \times E(\mathbf{x},\omega;\mathbf{x}_s) -\mathrm{i}\omega\mu(\mathbf{x})  H(\mathbf{x},\omega;\mathbf{x}_s)  = M(\mathbf{x}_s,\omega),\quad \mathbf{x}\in X, \omega\in \Omega \\
    \nabla \times H(\mathbf{x},\omega;\mathbf{x}_s)  -\sigma(\mathbf{x}) E(\mathbf{x},\omega;\mathbf{x}_s)  =J(\mathbf{x}_s,\omega) ,\quad  \mathbf{x}\in X, \omega\in \Omega
  \end{cases},
\end{equation}
where the variables $\mathbf{x}$, $t$ and $\omega$ denote space, time and frequency, respectively. The convention of Fourier transform $\partial_t \leftrightarrow -\mathrm{i}\omega$ has been adopted above. The electrical and magnetic fields ($E$ and $H$) are vectors consisting of 3 components in x, y and z directions; $J$ and $M$ stand for electrical and magnetic sources at the location $\mathbf{x}_s$, respectively. 
The magnetic permeability is $\mu$. The conductivity $\sigma$ is a symmetric $3\times 3$ 
tensor, i.e., $ \sigma_{ij}=\sigma_{ji}$, $i,j\in\{x,y,z\}$. In the isotropic medium, only the diagonal elements of the conductivity
tensor are non-zeros and the same in all directions:
$\sigma_{xx}=\sigma_{yy}=\sigma_{zz}$; $\sigma_{ij}=0, i\neq j $.
In the vertical transverse isotropic (VTI) medium, the diagonal 
elements are different in horizontal and vertical directions: $ \sigma_h:=\sigma_{xx}=\sigma_{yy},\quad \sigma_v=\sigma_{zz}$, where $\sigma_h$ and $\sigma_v$ stand for horizontal conductivity and vertical conductivity, respectively. The resistivity is defined as the inverse of the conductivity, i.e., $\rho_{ij}=1/\sigma_{ij}$.

The above diffusive Maxwell equation can be converted into wave domain \citep{Mittet_2010_HFD}, by defining a fictitious di-electrical permittivity in equation \eqref{eq:freq} as $\sigma = 2\omega_0 \varepsilon$,  while multiplying the 2nd equation in \eqref{eq:freq} with $\sqrt{-\mathrm{i}\omega/2\omega_0}$:
\begin{equation}\label{eq:mittet}
  \begin{cases}
\nabla\times E'(\mathbf{x},\omega;\mathbf{x}_s)-\mathrm{i}\omega'\mu H'(\mathbf{x},\omega;\mathbf{x}_s) = M'(\mathbf{x}_s,\omega),\quad \mathbf{x}\in X, \omega\in \Omega\\
\nabla\times H'(\mathbf{x},\omega;\mathbf{x}_s)+\mathrm{i}\omega'\varepsilon E'(\mathbf{x},\omega;\mathbf{x}_s)= J'(\mathbf{x}_s,\omega),\quad \mathbf{x}\in X, \omega\in \Omega
\end{cases} 
\end{equation}
based on the following correspondence relation
\begin{equation}
    E'=E,\;
    H'= \sqrt{\frac{-\mathrm{i}\omega}{2\omega_0}} H,\;
    M'= M,\;
    J'= \sqrt{\frac{-\mathrm{i}\omega}{2\omega_0}} J,\;
    \omega' = (1+\mathrm{i})\sqrt{\omega\omega_0}.
\end{equation}

The time domain counterpart of equation \eqref{eq:mittet} reads
\begin{equation}\label{eq:fictitiouswave}
  \begin{cases}
    \nabla\times E'(\mathbf{x},t';\mathbf{x}_s) + \mu(\mathbf{x})\partial_{t'} H'(\mathbf{x},t';\mathbf{x}_s) = M'(\mathbf{x}_s,t'), \quad  \mathbf{x}\in X, t'\in T\\
    \nabla\times H'(\mathbf{x},t';\mathbf{x}_s) - \varepsilon(\mathbf{x})\partial_{t'} E'(\mathbf{x},t';\mathbf{x}_s) = J'(\mathbf{x}_s,t'), \quad  \mathbf{x}\in X, t'\in T,
  \end{cases}
\end{equation}
which allows us to do efficient modelling using leap-frog finite-difference time-domain (FDTD) method over the staggered grid.  From the electromagnetic fields in the time domain, the frequency domain fields can be integrated during modelling using the  transformation 
\begin{equation}\label{eq:dtft}
  u(\mathbf{x},\omega; \mathbf{x}_s) = \int_0^{T_{\max}} u(\mathbf{x},t'; \mathbf{x}_s)e^{\mathrm{i}\omega't'} \mathrm{d}t' =  \int_0^{T_{\max}} u(\mathbf{x},t'; \mathbf{x}_s)e^{-\sqrt{\omega\omega_0}t'} e^{\mathrm{i}\sqrt{\omega\omega_0}t'} \mathrm{d}t',
\end{equation}
where $u(\mathbf{x},t';\mathbf{x}_s)\in \{E'_x,E'_y,E'_z, H'_x, H'_y, H'_z\}$ denotes the electric/magnetic field at the spatial location $\mathbf{x}\in X$ and the fictitious time $t'\in [0, T_{\max}]$ excited by a source from the spatial location $\mathbf{x}_s$; $T_{\max}$ is the final time that the frequency domain field $u(\mathbf{x},\omega;\mathbf{x}_s)$   reaches its steady state. Equation \eqref{eq:dtft} shows that the complex-valued frequency $\omega' = (1+\mathrm{i})\sqrt{\omega\omega_0}$ leads to an exponentially decay factor in the time integration, which is the key to realize the attenuation effect of diffusive EM field during time evolution.

To model the CSEM response as precise as possible, we resort to our newly developed high-order finite-difference time-domain method over non-uniform grid \citep{Yang_2023_HFDNU,Yang_2023_libEMM}.  To satisfy homogeneous boundary condition approximately, perfectly matched layers (PML) \citep{Komatitsch_2007_GEO} are padded surrounding the domain of interest to absorb the reflections in the truncated domain. The implementation of \eqref{eq:dtft} can be carried out on the fly during timestepping modelling thanks to the discrete time  Fourier transform (DTFT):
\begin{equation}\label{eq:dft}
  u(\mathbf{x},\omega_k;\mathbf{x}_s) = \sum_{n=0}^{N_t-1} u(\mathbf{x},t'_n;\mathbf{x}_s) \exp(-\sqrt{\omega_k\omega_0} t'_n) \exp(\mathrm{i}\sqrt{\omega_k\omega_0} t'_n),
\end{equation}
where the discretized fictitious time is $t'_n=n\Delta t$ with time step $\Delta t$ and time index $n=0,\cdots,N_t$, while the interested frequencies are $\omega_k$, $k=1,\cdots, N_{\omega}$. Note that the total number of time steps $N_t$ is much large than  the total number of discrete frequencies  $N_\omega$, i.e., $N_t\gg N_\omega$. The modelling terminates when the electric field $E'(\mathbf{x},\omega; \mathbf{x}_s)$ and the magnetic field $H'(\mathbf{x},\omega; \mathbf{x}_s)$ evolves to the stead state. This implies that we may check the convergence of the field regularly during the modelling to avoid additional timesteppings which have negligible contributions to the time integral. It is noteworthy that multiple frequencies can be integrated on the fly during the same timestepping procedure. These advantages make the method attractive for forward modelling.

\section{Fictitious wave domain CSEM inversion}

The CSEM inversion is an iterative nonlinear optimization procedure. At each iteration, the gradient of the data misfit with respect to conductivity can be computed via the product between forward and adjoint fields.
To use fictitious wave domain modelling engine for adjoint simulation, we present a novel approach to estimating the adjoint source time functions by formulating a regularized linear inverse problem. An in-depth analysis shows that this linear inverse problem needs to be solved only once to result in a limited number of basis functions, which can then be used to construct the time-domain adjoint source at different receiver locations.  This key development enables efficient and practical 3D CSEM inversion by fictitous wave domain modelling engine.

\subsection{CSEM inverse problem}

The Maxwell equations forms a linear system as
\begin{equation}\label{eq:matrixpde1}
  \underbrace{\begin{bmatrix}
      -\sigma &\nabla \times\\
      \nabla\times & -\mathrm{i}\omega\mu 
  \end{bmatrix}}_{\mathbf{A}(m)} \underbrace{\begin{bmatrix}
      E\\
      H
  \end{bmatrix}}_\mathbf{u}=\underbrace{\begin{bmatrix}
      J\\
      M
  \end{bmatrix}}_\mathbf{f},
\end{equation}
where the electromagnetic propagator $\mathbf{A}(m)$ is a linear operator applied to the vector field $\mathbf{u}$, which gathers electrical and magnetic fields (each component has  been denoted by $u:=u(\mathbf{x},\omega; \mathbf{x}_s)$ in the previous section). The  source vector with nonzero excitation at the source location $\mathbf{x}_s$ is prescribed by
$\mathbf{f}:=\mathbf{f}(\mathbf{x},\omega;\mathbf{x}_s)$. Note that all field variables are functions of both frequency and space ($u(\mathbf{x},\omega;\mathbf{x}_s), \mathbf{x}\in X, \omega\in \Omega$), while the medium property is only a function of space ($m(\mathbf{x}), \mathbf{x}\in X$), which can be the conductivity $\sigma_{ij}$ or the permeability $\mu$. 

To find the resistivity of the subsurface,  the data misfit is defined in least-squares sense
\begin{equation}\label{eq:phid}
  \phi_d(m)  =\frac{1}{2}\Big\|\underbrace{\begin{bmatrix}
      \mathbf{W}_1 &0\\
      0 & \mathbf{W}_2
  \end{bmatrix}}_\mathbf{W}\Big(\underbrace{\begin{bmatrix}
      E^{obs}(\mathbf{x}_r,\omega;\mathbf{x}_s)\\
      H^{obs}(\mathbf{x}_r,\omega;\mathbf{x}_s)
  \end{bmatrix}}_{\mathbf{d}(\mathbf{x}_r, \omega; \mathbf{x}_s)}-
  \underbrace{\begin{bmatrix}
      E(\mathbf{x}_r,\omega;\mathbf{x}_s)\\
      H(\mathbf{x}_r,\omega;\mathbf{x}_s)
  \end{bmatrix}}_{\mathbf{d}_{syn}(\mathbf{x}_r, \omega; \mathbf{x}_s)[m]} \Big) \Big\|^2
  = \frac{1}{2}\|\mathbf{W}(\mathbf{d}-\mathbf{R} \mathbf{u})\|^2,
\end{equation}
where  $\mathbf{d}(\mathbf{x}_r,\omega; \mathbf{x}_s)$ denotes the observed EM data  at the receiver location $\mathbf{x}_r$ due to the source at the location $\mathbf{x}_s$, while the synthetic data $\mathbf{d}_{syn}(\mathbf{x}_r,\omega;\mathbf{x}_s)[m]:=\mathbf{u}(\mathbf{x}_r,\omega;\mathbf{x}_s)=\mathbf{R}\mathbf{u}(\mathbf{x},\omega;\mathbf{x}_s)$ (simulated with the model parameter $m$) are extracted by the restriction operator $\mathbf{R}$ from the modelled wavefield at the receiver location
\begin{equation}
  u(\mathbf{x}_r,\omega; \mathbf{x}_s)=\int_X u(\mathbf{x},\omega;\mathbf{x}_s)\delta(\mathbf{x}-\mathbf{x}_r)\mathrm{d}\mathbf{x}, \quad \mathbf{x}\in X, \omega\in\Omega
\end{equation}
where $u(\mathbf{x}_r,\omega; \mathbf{x}_s)$ is one component of the vector field $\mathbf{u}(\mathbf{x}_r,\omega; \mathbf{x}_s)$.  In the remainder of the paper, we shall drop the dependence of $(\mathbf{x}_r,\omega;\mathbf{x}_s)$ and the summation over sources and receivers without loss of clarity unless clearly stated if necessary. Two weighting matrices, $\mathbf{W}_1$ and $\mathbf{W}_2$ together forming the diagonal weighting matrix $\mathbf{W}$, are employed to weight the electric and magnetic fields respectively. These weighting matrices may be specified according to the uncertainty model in terms of the real acquisition and the equipment \citep{mittet2012detection}. If the magnetic data are not considered in the inversion, we simply set $\mathbf{W}_2=0$.

Taking the first derivative of equation \eqref{eq:matrixpde1} with respect to model parameter $m$ gives
\begin{equation}
\frac{\partial \mathbf{u}}{\partial m} = -\mathbf{A}^{-1}(m)\frac{\partial \mathbf{A}(m)}{\partial m}\mathbf{u}.
\end{equation}
The gradient of the data misfit with respect to the model parameter $m$ is
\begin{equation}\label{eq:gradientderivation}
  \begin{split}
    \frac{\partial \phi_d(m)}{\partial m}=&\Re\langle \mathbf{W}\mathbf{R} \frac{\partial \mathbf{u}}{\partial m}, \mathbf{W}(\mathbf{R}\mathbf{u}-\mathbf{d})\rangle\\
    =&-\Re\langle \mathbf{W}\mathbf{R}\mathbf{A}^{-1}(m)\frac{\partial \mathbf{A}(m)}{\partial m}\mathbf{u}, \mathbf{W}(\mathbf{R}\mathbf{u}-\mathbf{d})\rangle\\
    =&\Re\langle \frac{\partial \mathbf{A}(m)}{\partial m}\mathbf{u},
    \underbrace{(\mathbf{A}^\dagger)^{-1}(m)\mathbf{R}^\dagger\mathbf{W}^\dagger\mathbf{W}(\mathbf{d}-\mathbf{R}\mathbf{u}}_{\mathbf{v}})\rangle\\
    =&\Re\langle \frac{\partial \mathbf{A}(m)}{\partial m}\mathbf{u},\mathbf{v}\rangle\\
    =&\Re \sum_s \sum_\omega \bar{\mathbf{v}}^\mathrm{T}\frac{\partial \mathbf{A}(m)}{\partial m} \mathbf{u},
  \end{split}
\end{equation}
where $\Re$ takes the real part of a complex number, $\dagger$ is the complex conjugate transpose. It is noteworthy that the newly introduced variable $\mathbf{v}$, coined adjoint variable, or co-state variable in optimal control theory, must satisfy
\begin{equation}\label{eq:adj0}
  \mathbf{A}^\dagger(m) \mathbf{v}=\mathbf{R}^\mathrm{T} \mathbf{W}^\mathrm{T} \mathbf{W}(\mathbf{d}-\mathbf{R}\mathbf{u}).
\end{equation}
This shows that $\mathbf{v}$ satisfies another Maxwell equation based on adjoint Maxwell operator $\mathbf{A}^\dagger(m)$. The right hand side of the adjoint equation acts as the virtual source to emanate the adjoint field. 

Equation \eqref{eq:gradientderivation} is a generic expression applicable in fully anisotropic medium for the model parameters $m\in\{\sigma_{ij},\mu\}$. Let us point out that the above gradient expression may be scaled with the local cell volume $\Delta V(\mathbf{x})$ due to the discretization of the spatial integral $\int_X \mathrm{d}\mathbf{x}$.
In the CSEM settings, the magnetic permeability $\mu$ is considered as the same constant as in the vacuum. We are interested in retrieving the conductivity $\sigma$ (or the resistivity $\rho=1/\sigma$) which may exhibit anisotropy.  Equation \eqref{eq:matrixpde1} implies
\begin{equation}
\frac{\partial \mathbf{A}(m)}{\partial \sigma}=\frac{\partial}{\partial \sigma}\begin{bmatrix}
      -\sigma &\nabla \times\\
      \nabla\times & -\mathrm{i}\omega\mu 
  \end{bmatrix}=\begin{bmatrix}
      -1 &0\\
      0 & 0 
  \end{bmatrix}.
\end{equation}
Denote the adjoint field $\mathbf{v}=(\underline{E},\underline{H})^\mathrm{T}$, where the underline is utilized to distinguish it from the forward field. The gradient  in equation \eqref{eq:gradientderivation} becomes
\begin{equation}\label{eq:gradsigmaij}
  \frac{\partial \phi_d}{\partial\sigma_{ij}} = -  \Re\sum_s \sum_\omega \underline{\bar{E}}_j \cdot E_i, \qquad i,j\in\{x,y,z\}.
\end{equation}
For commonly used VTI medium, the gradients of the misfit with respect to horizontal conductivity $\sigma_h$ and vertical conductivity $\sigma_v$ are
\begin{equation}
  \frac{\partial\phi_d}{\partial \sigma_h} =-\Re \sum_s \sum_\omega( \underline{\bar{E}}_x \cdot E_x + \underline{\bar{E}}_y \cdot E_y), \quad
  \frac{\partial\phi_d}{\partial \sigma_v} =-\Re \sum_s \sum_\omega \underline{\bar{E}}_z \cdot E_z.
\end{equation}
In isotropic medium ($\sigma_{ii}=\sigma$, $i=x,y,z$), we have the simple gradient expression
\begin{equation}
  \frac{\partial \phi_d}{\partial \sigma} =-\Re \sum_s \sum_\omega \underline{E}^\dagger\cdot E
  =-\Re\sum_s \sum_\omega (\underline{\bar{E}}_x\cdot E_x + \underline{\bar{E}}_y\cdot E_y+\underline{\bar{E}}_z\cdot E_z).
\end{equation}

\subsection{Adjoint modelling in fictitious wave domain}

To facilitate the computation of the adjoint field and the inversion gradient, one may take the conjugate of equation \eqref{eq:adj0}:
\begin{equation}\label{eq:adj1}
  \mathbf{A}^\mathrm{T}(m) \bar{\mathbf{v}}=\overline{\mathbf{R}^\mathrm{T} \mathbf{W}^\mathrm{T} \mathbf{W}(\mathbf{d}-\mathbf{R}\mathbf{u})}.
\end{equation}
Since $\mathbf{A}^\mathrm{T}(m)=\mathbf{A}(m)$, the above operation allows us to use the same modelling engine to compute the conjugate of the adjoint field, as long as proper adjoint source can be provided for adjoint simulation.
Unfortunately, for CSEM inversion in the frequency domain, only few discrete frequencies are available, leading to the adjoint source (i.e. the data residual) specified only at few discrete frequencies. A time domain adjoint source is therefore not directly accessible.

To still use fictitious wave domain timestepping modelling for adjoint simulation, let us now split the weighted data residual into electrical and magnetic components $\overline{\mathbf{R}^\mathrm{T} \mathbf{W}^\mathrm{T} \mathbf{W}(\mathbf{d}-\mathbf{R}\mathbf{u})}=(\overline{\delta d_E}, \overline{\delta d_H})^\mathrm{T}$ and repeat the substitution and multiplication operations as in Section \ref{sec:forward}:
\begin{equation}\label{eq:adj2}
  \begin{cases}
    \nabla\times \underbrace{\overline{\underline{E}}}_{\underline{\overline{E}}'} +\mu\underbrace{\sqrt{-\mathrm{i}2\omega \omega_0}}_{-\mathrm{i}\omega'}\underbrace{\sqrt{\frac{-\mathrm{i}\omega}{2\omega_0}} \underline{\overline{H}}}_{\underline{\overline{H}}'} =\overline{ \delta d_E},\\
    \nabla\times\underbrace{ \sqrt{\frac{-\mathrm{i}\omega}{2\omega_0}} \underline{\overline{H}}}_{\underline{\overline{H}}'} -\underbrace{\sqrt{-\mathrm{i}2\omega\omega_0}}_{-\mathrm{i}\omega'}\varepsilon  \underline{\overline{E}} = \sqrt{\frac{-\mathrm{i}\omega}{2\omega_0}} \overline{\delta d_H}.
    \end{cases}
\end{equation}
We need to consider the time domain counterpart of the above system for efficient adjoint modelling by fictitious timestepping modelling. Switching from diffusive frequency domain back to fictitious time domain with wave should be understood as a linear  transformation rather than the true inverse Fourier transform. This is because the mapping of fictitious transform is not bijective, as the number of time steps $N_t$ for numerical simulation is significantly larger than the number of discrete frequencies $N_\omega$ used for CSEM investigation ($N_\omega\ll N_t$). There exists an infinite number of time series which may match the few discrete frequencies. Searching for the inverse is therefore an under-determined problem.

Let us denote the right hand side of equation \eqref{eq:adj2} ($\overline{\delta d_E}$ or $\sqrt{\frac{-\mathrm{i}\omega}{2\omega_0}}\overline{\delta d_H}$) at a specific receiver location as $s(\omega)$. To retrieve a long time series for adjoint modelling, a linear inverse problem is formulated to convert the adjoint source $s(\omega)$  from frequency to time domain based on DTFT in equation \eqref{eq:dft} by minimizing the following misfit functional
\begin{equation}\label{eq:leastsquareadjsrc}
  \begin{split}
    \Psi &= \sum_{k=1}^{N_\omega} |s(\omega_k) -\sum_{n=0}^{N_t-1} e^{-\sqrt{\omega_k\omega_0} t'_n} e^{\mathrm{i}\sqrt{\omega_k\omega_0} t'_n} s(t'_n)|^2 + \gamma \sum_{n=0}^{N_t-1}|s(t'_n)|^2 \nonumber \\
    & =\sum_{k=1}^{N_\omega}\Big (|\Re\{ s(\omega_k)\} -\sum_{n=0}^{N_t-1} e^{-\sqrt{\omega_k\omega_0} t'_n}\cos(\sqrt{\omega_k\omega_0} t'_n) s(t'_n)|^2 \nonumber \\
    &\qquad + |\Im\{ s(\omega_k)\} -\sum_{n=0}^{N_t-1} e^{-\sqrt{\omega_k\omega_0} t'_n}\sin(\sqrt{\omega_k\omega_0} t'_n) s(t'_n)|^2 \Big)+\gamma \sum_{n=0}^{N_t-1}|s(t'_n)|^2
  \end{split}
\end{equation}
where  $\Im$ takes the imaginary part of the complex variable. If we  discretize the time as $t_n=n\Delta t$ and define
\begin{equation}
  \mathbf{s}(\omega):=\begin{bmatrix}
  \Re\{s(\omega_1)\}\\
  \vdots\\
  \Re\{s(\omega_{N_\omega})\}\\
  \Im\{s(\omega_1)\}\\
  \vdots\\
  \Im\{s(\omega_{N_\omega})\}
  \end{bmatrix}, \mathbf{s}(t'):= \begin{bmatrix}
    s(t'_0)\\
    \vdots\\
    s(t'_{N_t-1})
  \end{bmatrix}
\end{equation}
and 
\begin{equation}
  \mathbf{B}:=\begin{bmatrix}
  e^{-\sqrt{\omega_1\omega_0} t'_0}\cos(\sqrt{\omega_1\omega_0} t'_0) & \cdots & e^{-\sqrt{\omega_1\omega_0} t'_{N_t-1}}\cos(\sqrt{\omega_1\omega_0} t'_{N_t-1})\\
  \vdots & \ddots & \vdots \\
  e^{-\sqrt{\omega_{N_\omega}\omega_0} t'_0}\cos(\sqrt{\omega_{N_\omega}\omega_0} t'_0) & \cdots & e^{-\sqrt{\omega_{N_\omega}\omega_0} t'_{N_t-1}}\cos(\sqrt{\omega_{N_\omega}\omega_0} t'_{N_t-1})\\
  e^{-\sqrt{\omega_1\omega_0} t'_0}\sin(\sqrt{\omega_1\omega_0} t'_0) & \cdots & e^{-\sqrt{\omega_1\omega_0} t'_{N_t-1}}\sin(\sqrt{\omega_1\omega_0} t'_{N_t-1})\\
  \vdots & \ddots & \vdots \\
  e^{-\sqrt{\omega_{N_\omega}\omega_0} t'_0}\sin(\sqrt{\omega_{N_\omega}\omega_0} t'_0) & \cdots & e^{-\sqrt{\omega_{N_\omega}\omega_0} t'_{N_t-1}}\sin(\sqrt{\omega_{N_\omega}\omega_0} t'_{N_t-1})\\
  \end{bmatrix},
\end{equation}
the objective of the above linear inverse problem can be compactly written in matrix form
\begin{equation}
  \Psi = \|\mathbf{s}(\omega)-\mathbf{B} \mathbf{s}(t')\|^2 +\gamma\|\mathbf{s}(t')\|^2,
\end{equation}
where  a regularization term taking into account the minimum energy of the solution has been penalized by a parameter $\gamma$ to stabilize  the linear inversion and to determine a unique solution.  Note that the complex factors have been split into real and imaginary part. To solve this ill-posed problem, regularization has been introduced to recover the well-posedness during the inversion of the matrix $\mathbf{B}\in \mathbb{R}^{2N_\omega \times N_t}$.
The solution $\mathbf{s}(t)$ can easily be found using linear optimization algorithms, e.g. LSQR \citep{Paige_1982_ALS} and the conjugate gradient method for the resulting normal equation (CGNR) \citep{Saad_2003_IMS}. Based on the frequency and time index, the each element of $\mathbf{B}$ can be formed on the fly when computing matrix vector product, implying a matrix free implementation.

It should be noted that fictitious wave domain modelling is simply a mathematical tool to efficiently compute frequency domain EM fields. Disguising a diffusive phenomenon as a wave event  will not change the diffusive nature of the underlying physics. Indeed, \citet{Mittet_2010_HFD} proposed another method to calculate the adjoint source time function by superposition of delayed causal wavelets. However, it yields a nonlinear inverse problem which is very difficult to solve. Clearly, embedding a linear inverse problem in a nonlinear inverse problem is a much better option than embedding double nonlinear inverse problems. The truncated singular value decomposition (SVD) have been applied \citep{Storen_2008_Gradient} to invert  for an adjoint source, following the modified wave domain formulation of \citet{Maao_2007_FFT}.
Mathematically, our approach using damped least-squares should give the same solution as the one from truncated SVD, since SVD is equivalent to finding the Moore-Penrose pseudo-inverse for $\mathbf{B}$  \citep[theorem 1.2.10 in p. 15, and section 2.7.2, p. 101]{Bjorck_1996_NML}. However, our empirical experience shows that the approximate inverse given by SVD suffers from numerical instability while consuming too much memory.  The above linear inversion method gets rid of these issues,  incorporating all frequencies in one run with better memory and computational efficiency.

\subsection{A matrix-free basis function implementation}

Let us denote $\mathbf{B}^+=(\mathbf{B}^\mathrm{T} \mathbf{B} + \gamma \mathbf{I})^{-1} \mathbf{B}^\mathrm{T} \in \mathbb{R}^{N_t\times 2N_\omega}$ the Moore-Penrose pseudo-inverse of $\mathbf{B}\in \mathbb{R}^{2N_\omega\times N_t}$ such that
\begin{equation}\label{eq:mpinv}
  \mathbf{s}(t) = \mathbf{B}^+\mathbf{s}(\omega).
\end{equation}
The matrix $ \mathbf{B}^+$ may be written down using the $2N_\omega$ columns:
\begin{equation}
  \mathbf{B}^+  =\begin{bmatrix}
  \mathbf{b}_1,\cdots,\mathbf{b}_{N_\omega},\mathbf{b}_{N_\omega+1},\cdots,\mathbf{b}_{2N_\omega}
  \end{bmatrix}
\end{equation}
in which  $\mathbf{b}_k\in \mathbb{R}^{N_t}$ is the $i$-th column of $\mathbf{B}^+$. As a result, equation \eqref{eq:mpinv} translates into
\begin{equation}\label{eq:basis}
  \mathbf{s}(t)=\sum_{k=1}^{N_\omega} \Re\{ s(\omega_k)\} \mathbf{b}_k  + \Im\{ s(\omega_k) \} \mathbf{b}_{k+N_\omega}.
\end{equation}
It then becomes clear that the columns of $\mathbf{B}^+$  are the very few number of basis functions to construct the adjoint source time function. Given $\Delta t$ and $N_t$, these basis functions are uniquely determined, hence independent of the frequency spectrum $s(\omega)$ (which are simply the coefficients for the linear combination of these basis functions).
This means the iterative solution procedure needs to be performed only once to find all $\mathbf{b}_k$, which can then be used for each receiver locations to infer the corresponding adjoint source time function. Equation \eqref{eq:basis} also gives a recipe to efficiently compute these functions:
\begin{itemize}
\item If we set
  $\Re\{s(\omega_k)\}=1$, $\Im\{s(\omega_k)\}=0$ and $\Re\{s(\omega_j)\}=\Im\{s(\omega_j)\}=0$ ($j\neq k$), we obtain $\mathbf{s}(t)=\mathbf{b}_k$.

\item If we set
  $\Im\{s(\omega_k)\}=1$, $\Re\{s(\omega_k)\}=0$ and $\Re\{s(\omega_j)\}=\Im\{s(\omega_j)\}=0$ ($j\neq k$), we obtain $\mathbf{s}(t)=\mathbf{b}_{k+N_\omega}$.
\end{itemize}
By repeatedly feeding the iterative optimization algorithm $2N_\omega$ times using different Dirac delta input $\mathbf{s}(\omega)=(0,\cdots,0,1,0,\cdots, 0)^\mathrm{T}$, we obtain the $2N_{\omega}$ columns of $\mathbf{B}^+$.

\begin{figure}
  \centering
  \includegraphics[width=0.75\linewidth]{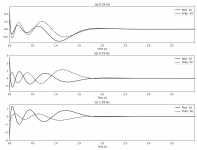}
  \caption{The basis functions using the parameters  $\Delta t=0.004$ s, $N_t=1000$ and $N_\omega=3$. The amplitude of these basis functions decays to zero with the increase of the time, which is quite suitable for frequency domain EM fields converging to steady state after simulation using sufficient number of time steps.}\label{fig:basisfunction}
\end{figure}

\begin{figure}
  \centering
  \includegraphics[width=0.75\linewidth]{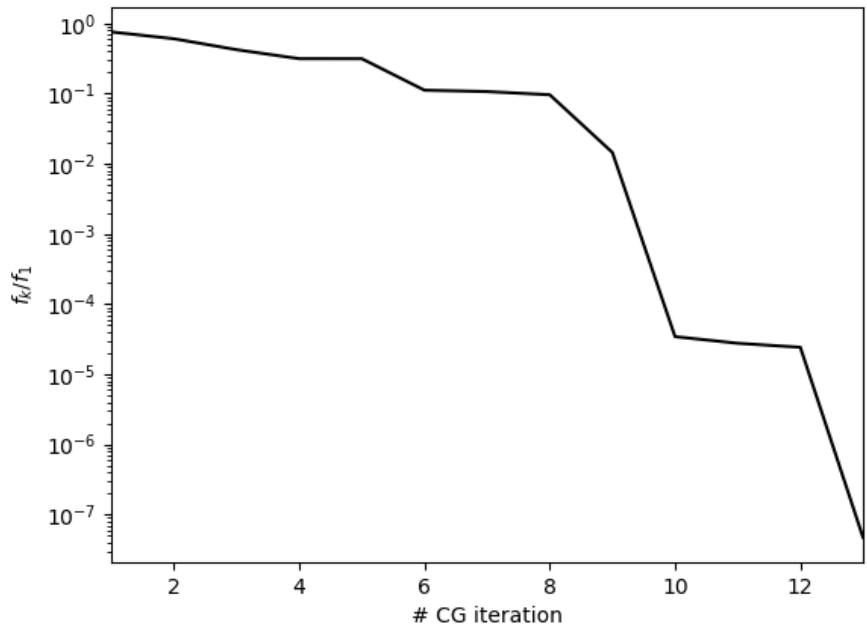}
  \caption{The convergence of CG iterations for the solution of basis functions}\label{fig:cgconv}
\end{figure}

Figure~\ref{fig:basisfunction} gives an example of these basis functions using the parameters $\Delta t=0.002$, $N_t=1000$ and $f_0=1$ and $N_\omega=3$. To solve this linear optimization problem, Claerbout's conjugate gradient algorithm \citep[chapter 2.3.6]{claerbout2008image} (which is essentially a CGNR algorithm according to \citet[chapter 8.3, pp. 266-268]{Saad_2003_IMS}) has been applied.
As can be seen from Figure~\ref{fig:cgconv}, within 15 CG iterations, the proposed method leads to highly accurate estimation for these basis functions, with the error less than $10^{-8}$. In this test, the regularization parameter $\gamma$ was chosen to be $10^{-3}$. It is found the shape of the basis functions do not change much using different choices of $\gamma$, indicating that the method is robust.

For different receiver locations, the data residuals in the frequency domain may be dramatically different, but these basis functions are the same. We therefore only solve the above regularized linear inverse problem once to construct all adjoint source time functions through their linear combination.
The frequency domain data residuals are the weighting coefficients to modulate these basis.

One thing to remark is that the adjoint modelling in the time domain should have reverse time order compared with the forward modelling. Taking the conjugate reverse the time again. This is the reason why a decaying trend of the estimated basis function is observed in Figure~\ref{fig:basisfunction}.

\subsection{Connections and distinctions with existing methods}

Note the incident field $E_i$ is the product between the Green's function with the source current
\begin{equation}
    E_i(\mathbf{x},\omega; \mathbf{x}_s)=\sum_k G_{ik}^{EE}(\mathbf{x},\omega; \mathbf{x}_s) J_k(\mathbf{x}_s,\omega),
\end{equation}
while the conjugate of adjoint field can be specified as
\begin{equation}\label{eq:conjadj}
    \overline{\underline{E}}_j(\mathbf{x},\omega; \mathbf{x}_s) =\int_X (\sum_k G_{jk}^{EE}(\mathbf{x},\omega; \mathbf{x}_r) \overline{\delta d^E_k(\mathbf{x}_r,\omega;\mathbf{x}_s)}+\sum_i G_{ji}^{EH}(\mathbf{x},\omega; \mathbf{x}_r) \overline{\delta d^H_i(\mathbf{x}_r,\omega;\mathbf{x}_s)})\mathrm{d}\mathbf{x}_r, 
\end{equation}
where $i, j, k\in\{x,y,z\}$ indicate different components of the EM field, while $G_{ij}^{EE}(\mathbf{x},\omega;\mathbf{x}_s)$ and $G_{ik}^{EH}(\mathbf{x},\omega;\mathbf{x}_s)$ stand for the $i$th component of the electrical Green's function for angular frequency $\omega$ at the spatial location $\mathbf{x}$ due to the electrical and magnetic sources directed in $j$- and $k$- th directions.
This implies that we can still use fictitious wave modelling for adjoint simulation: one can, for each frequency and each receiver, first compute the Green's function (which is independent of the source time function) in frequency domain, and then form the adjoint field by linear combination of them using the data residual. This method, ensured by the superposition principle thanks to the linearity of the Maxwell equation, enables the approach still working out. Such a computing scheme repeats many times of simulation depending on the number of receivers, which is very inefficient when the number of receivers becomes extremely large.

Equation \eqref{eq:adj1} shows that the adjoint field can be computed using only one modelling, provided that a time-domain adjoint source time function is available and then injected at once to do time-domain simulation. The spirit of our method is to perform frequency domain inversion using time domain modelling engine, similar to  \citet{Sirgue_2008_FDW,sirgue2010system} for seismic full waveform inversion (FWI). However, there are a number of differences. Switching between time and frequency domain is straightforward in seismic FWI using the definition of discrete inverse Fourier transform. One can convert frequency domain data residual using the definition of inverse discrete Fourier transform
\begin{equation}\label{eq:sirgue1}
  s(t) = \frac{1}{2\pi}\int_{-\infty}^\infty s(\omega)e^{-i\omega t} \mathrm{d}\omega \approx \sum_{k=1}^{N_\omega} s(\omega_k) e^{-i\omega_k t}.
\end{equation}
Of course, the same issue presented in \citet{sirgue2010system} approach: only  limited number of frequencies are available to reconstruct a long time series.
The last equality in \eqref{eq:sirgue1} assumes all the absent frequencies are zeros. Since Fourier basis is orthonormal, the computation of adjoint source time function under this assumption leads to a minimum energy solution  equivalent to solving the following least-squares minimization problem
\begin{equation}
  \min_{s(t)} \sum_{k=1}^{N_\omega}\|s(\omega_k)-\sum_{j=1}^{N_t} s(t_j) e^{\mathrm{i}\omega_k t_j}\|^2.
\end{equation}
For fictitious wave domain approach, plugging the complex-valued frequency $\omega'=(1+\mathrm{i})\sqrt{\omega\omega_0}$ into exponential factor yields an exponentially decay factor 
\begin{equation}
  e^{\mathrm{i}\omega't'}=e^{\mathrm{i}(1+\mathrm{i})\sqrt{\omega\omega_0}t'}=e^{-\sqrt{\omega\omega_0}t'}e^{\mathrm{i}\sqrt{\omega\omega_0}t'},
\end{equation}
leaving the converted EM field possessing strong attenuation behavior, while making the basis no more orthogonal. Precisely speaking, the transformation in \eqref{eq:dtft} is a Laplace transform rather than normal Fourier transform. As a result, the inverse fictitious transformation is not well defined, that is,
\begin{equation}\label{eq:fact}
 u(\mathbf{x},t'; \mathbf{x}_s) \neq \frac{1}{2\pi}\int_\Omega \mathrm{d}\omega u(\mathbf{x},\omega; \mathbf{x}_s) e^{-\mathrm{i}\omega't'}.
\end{equation}
 Since direct use of inverse Fourier transform as done by  \citet{sirgue2010system} does not apply, a matrix-free iterative solution of adjoint source time function in this work is evidently desirable.

\subsection{The final algorithm}

The total misfit functional for CSEM inversion consists of both data fitting and model regularization:
\begin{equation}
  \phi(m) = \phi_d(m) + \beta \phi_m(m),
\end{equation}
where $\phi_m(m)$, penalized by the parameter $\beta$, is the model misfit term to enforce smoothness on the inverted model. A popular choice is to use Tikhonov regularization to minimize the roughness of the model compared with a reference model $m_{ref}$, i.e., $  \phi_m(m) = \frac{1}{2}\|\nabla^\alpha (m-m_{ref})\|^2$,
where the anisotropic first order difference operator $\nabla^\alpha=(\alpha_x\partial_x, \alpha_y\partial_y,\alpha_z\partial_z)^\mathrm{T}$ with coefficients $\alpha_x$, $\alpha_y$ and $\alpha_z$ can be tuned by the user. This leads to the total gradient of the misfit functional prescribed by
$  \frac{\partial \phi}{\partial m} =\frac{\partial \phi_d}{\partial m} +\frac{\partial \phi_m}{\partial m}$.
Based on the gradient information computed above, a descent direction $\delta m^k$ can be constructed to update the model parameters iteratively
\begin{equation}
  m^{k+1} = m^k + \alpha \delta m^k,
\end{equation}
where  $\alpha$ is the step length estimated by line search method. In this paper, the descent direction $\delta m^k$ at the $k$-th iteration is computed using l-BFGS algorithm \citep{Nocedal_2006_NOO} by storing the gradients in the previous iterations.

It is worth noting that the inversion gradient computed using the proposed approach is approximate rather than exact. Consequently, the descent direction estimated using l-BFGS after several iterations may have difficulty to succeed in the line search procedure. We therefore propose to restart l-BFGS using the steepest descent direction after the first failure of the line search. The algorithm will eventually be terminated if the maximum number of iterations is reached.

Using conductivity as model parameter creates some kinds of ill-conditioning due to its value varying in a large range (the seawater is around 0.3 $\Omega\cdot$m while a thin layer of hydrocarbon bearing sediment is more than 100 $\Omega\cdot$m). To capture the high contrast of the parameter variations, the inversion is re-parameterized using the logarithmic transformation $m := \ln\rho$,
which is a dimensionless parameter with comparable magnitude varying in a much smaller dynamic range. Indeed, there exists other types of parameter scaling, see another example in \citet{abubakar20082}.
This re-parametrization gives a relation $\rho=e^m$. Switching the parametrization is trivial thanks to the chain rule:
\begin{equation}
  \frac{\partial\phi}{\partial m} = \frac{\partial \rho}{\partial m}\cdot\frac{\partial \sigma}{\partial \rho}\cdot \frac{\partial\phi}{\partial\sigma}=-\frac{1}{\rho} \frac{\partial\phi}{\partial\sigma}.
\end{equation}
During the iterative inversion,  the model parameter may be bounded within an upper bound $m_{\max}$ and a lower bound $m_{\min}$ according to a priori knowledge. These bounds help to stabilize the inversion in a more physically sensible manner.

\section{Application examples}

In this section we present two numerical examples to demonstrate the application of our method. The two examples are designed to validate the gradient expressions for the isotropic (the first example) and the VTI anisotropic (the second example) cases. In both examples, I use in-line and azimuth data with $E_x$ and $E_y$ components. This motivates us to utilize the uncertainty model proposed in \citet{morten2009uncertainy} to compute the diagonal weight matrix $\mathbf{W}$ in order to capture the varying sensitivity of broadside EM data. The observed data are poluted by 3\% of Gaussian white noise before inversion. The regularization parameters are configured with $\beta=0.01$ using a cooling factor (between 0.7 and 1) through iterations. We set $\alpha_x=\alpha_y=1$ and $\alpha_z=0.1$ in both tests. We choose the initial model as the reference model for both inverse exercises.  To mitigate the weak sensitivity of the EM field to the deeper part of the model, a depth preconditioning has been applied following the work of \citet{Plessix_2008_RIC}.
The code is parallelized over the sources  to achieve the best scaling performance since the sources are independent of each other.

\subsection{A two-block model}

\begin{figure}
  \centering
  \includegraphics[width=0.75\linewidth]{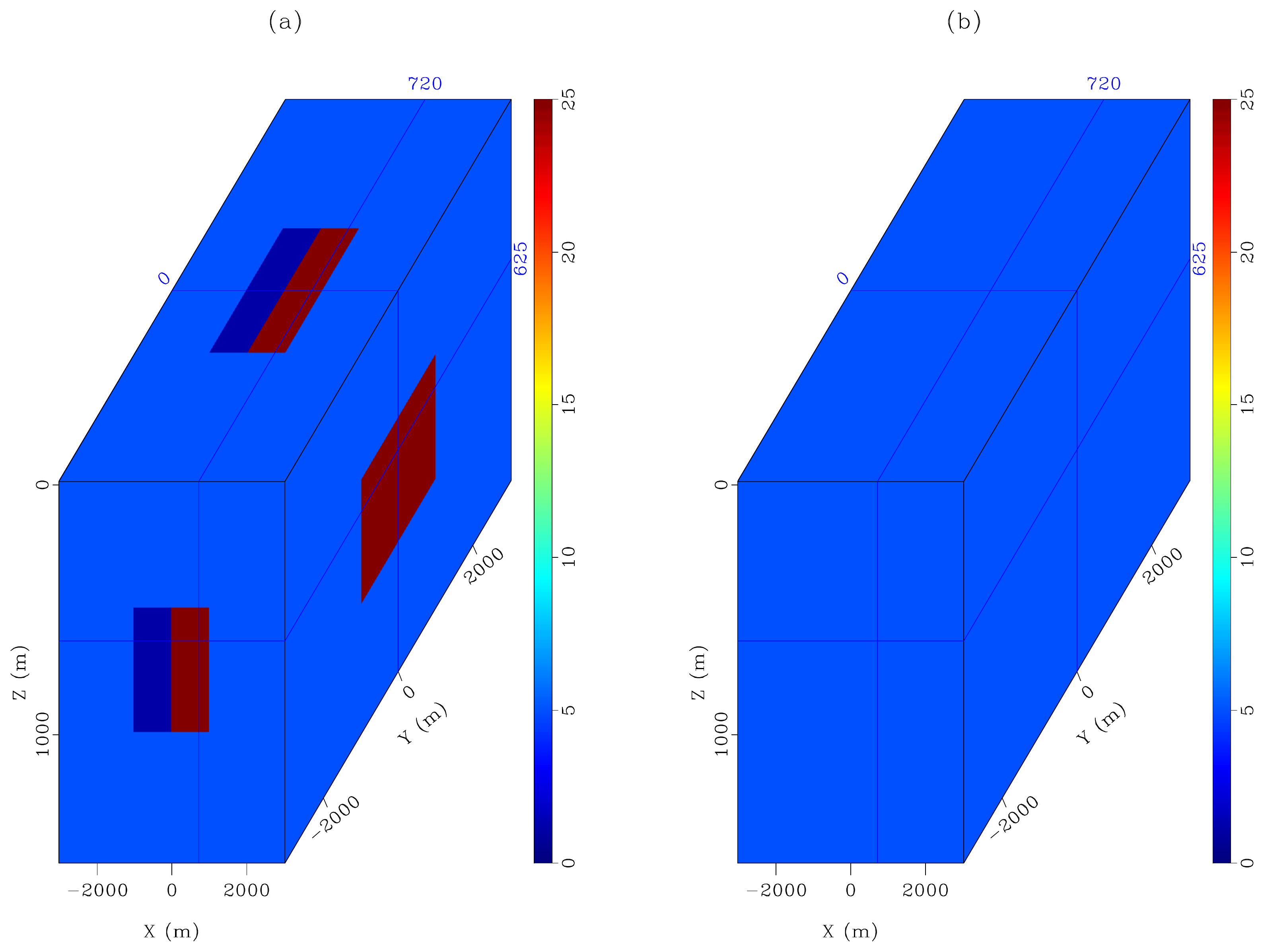}
  \caption{(a) The true land model (5 $\Omega\cdot$m background with two anomalous inclusions of 1 $\Omega\cdot$m (left) and 100 $\Omega\cdot$m (right) ); (b) The homogeneous initial model of 5 $\Omega\cdot$m. The models are clipped at 25 $\Omega\cdot$m for display purpose.}\label{fig:landmodel}
\end{figure}

\begin{figure}
  \centering
  \includegraphics[width=0.7\linewidth,angle=-90]{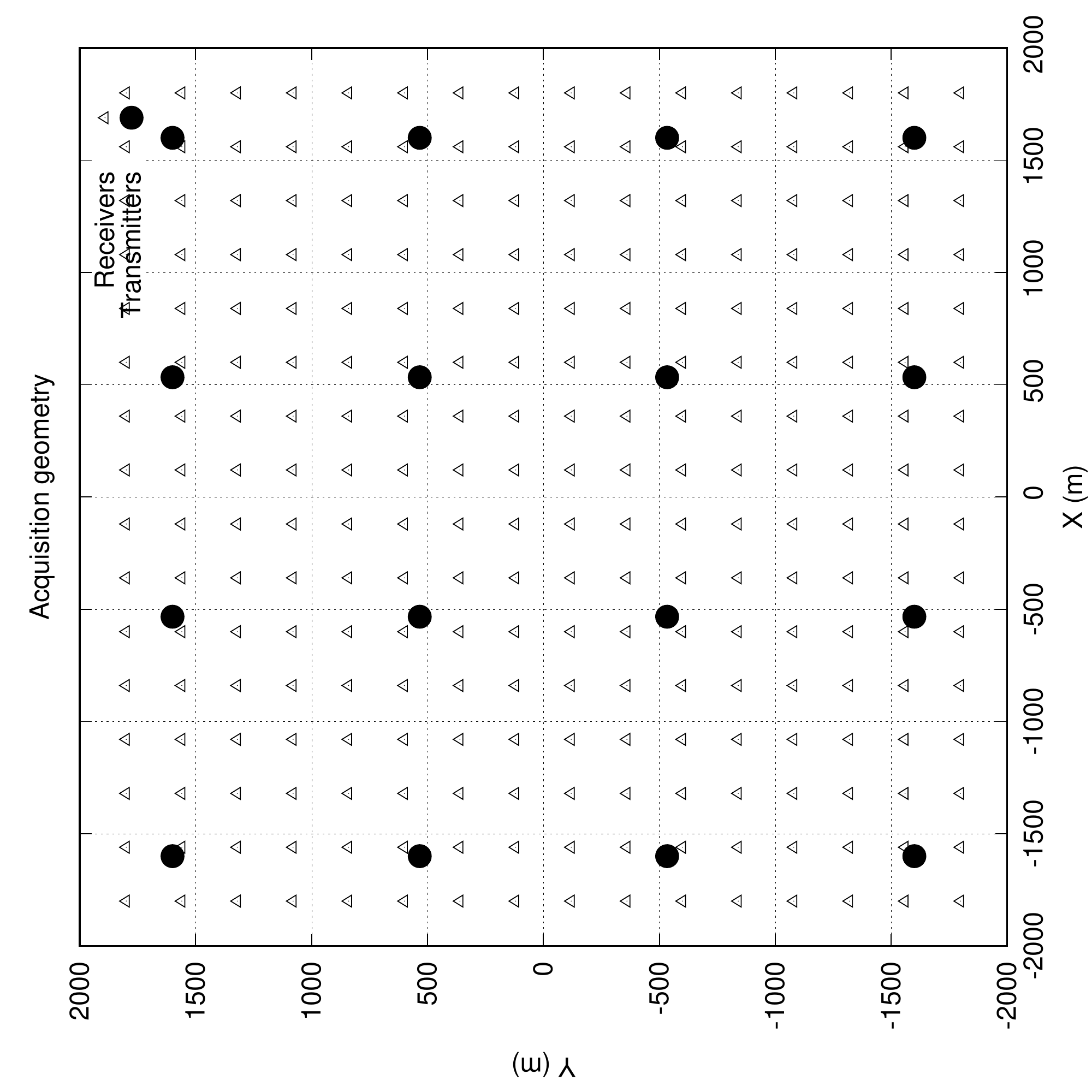}
  \caption{Survey layout sheet for the two-block model. In total, 16 transmitters (marked by dots) and 256 receivers (labeled with triangles) are deployed. }\label{fig:landsurvey}
\end{figure}

Our first example is a two-block land model similar to the one presented in \citet{grayver2013gji}. The model has a background resistivity of 5 $\Omega\cdot$m, including two anomalies: a low resistivity inclusion of 1 $\Omega\cdot$m and a high resistivity inclusion of 100 $\Omega\cdot$m, as shown in Figure~\ref{fig:landmodel}a. The physical dimension of the model spans over a 3D domain $X=X_1\times X_2\times X_3$, where $X_1=[-2000,2000]$ m, $X_2= [-2000,2000]$ m and $X_3= [0, 1500]$ m. There are 16 receivers and 256 transmitters with equal distance on the surface of this model. Figure~\ref{fig:landsurvey} gives the survey layout sheet of the acquisition geometry. The model has been discretized with the regular grid using $\Delta x=\Delta y=60$ m and $\Delta z=25$ m. The observed data are generated at two frequencies: 0.25 Hz and 1 Hz.

\begin{figure}
  \centering
  \includegraphics[width=0.8\linewidth]{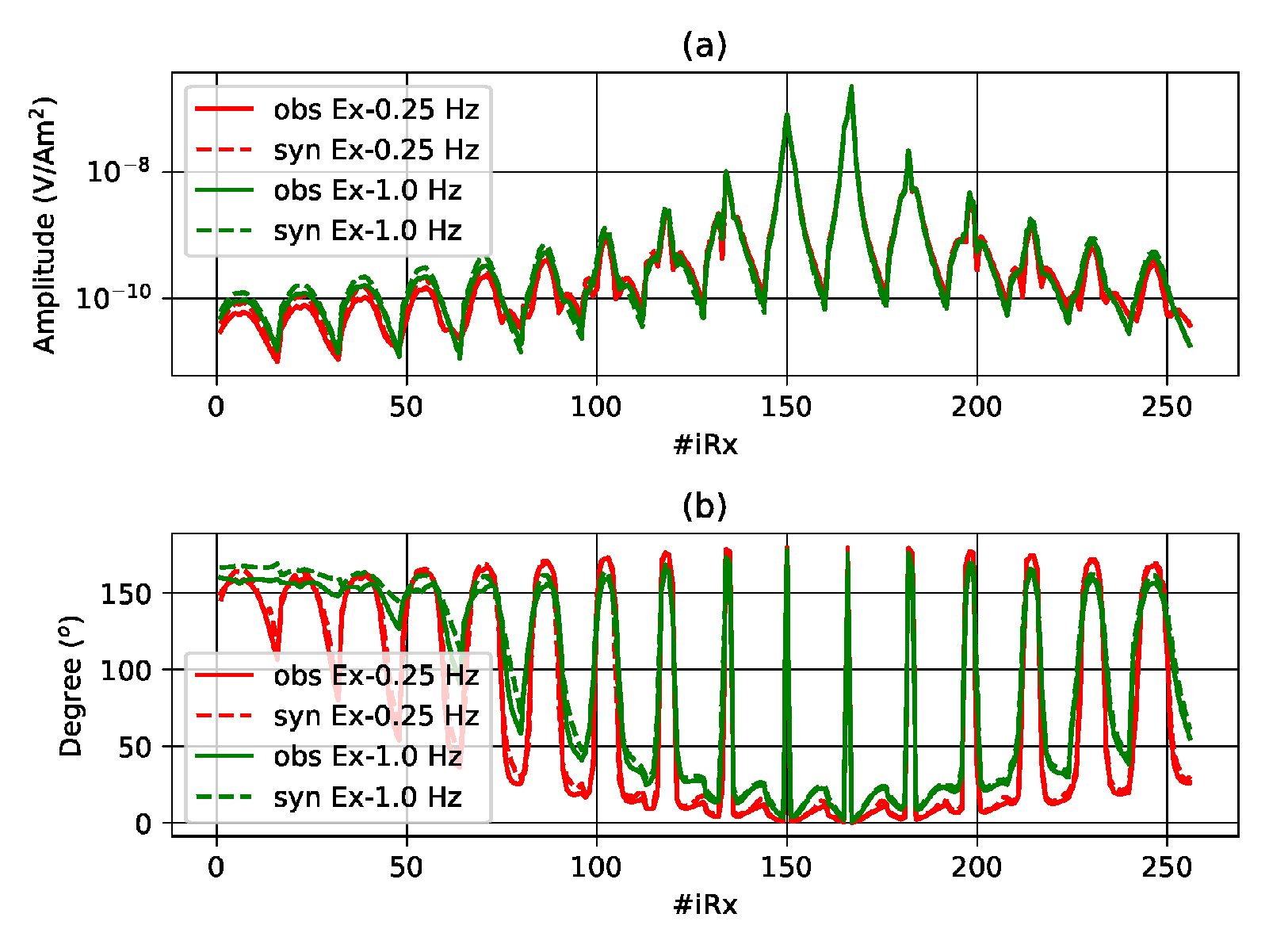}
  \caption{Comparison between observed data and synthetic data from initial model according to (a) amplitude and (b) phase for transmitter Tx-10.  The horizontal axis is the index of the receivers rather than offset, since we consider broadside configuration. It is clear that both the amplitude and phase does not match very well. Note also that the level of the data matching has a strong correlation with the receiver indices.}\label{fig:dataiter1}
\end{figure}

\begin{figure}
  \centering
  \includegraphics[width=0.8\linewidth]{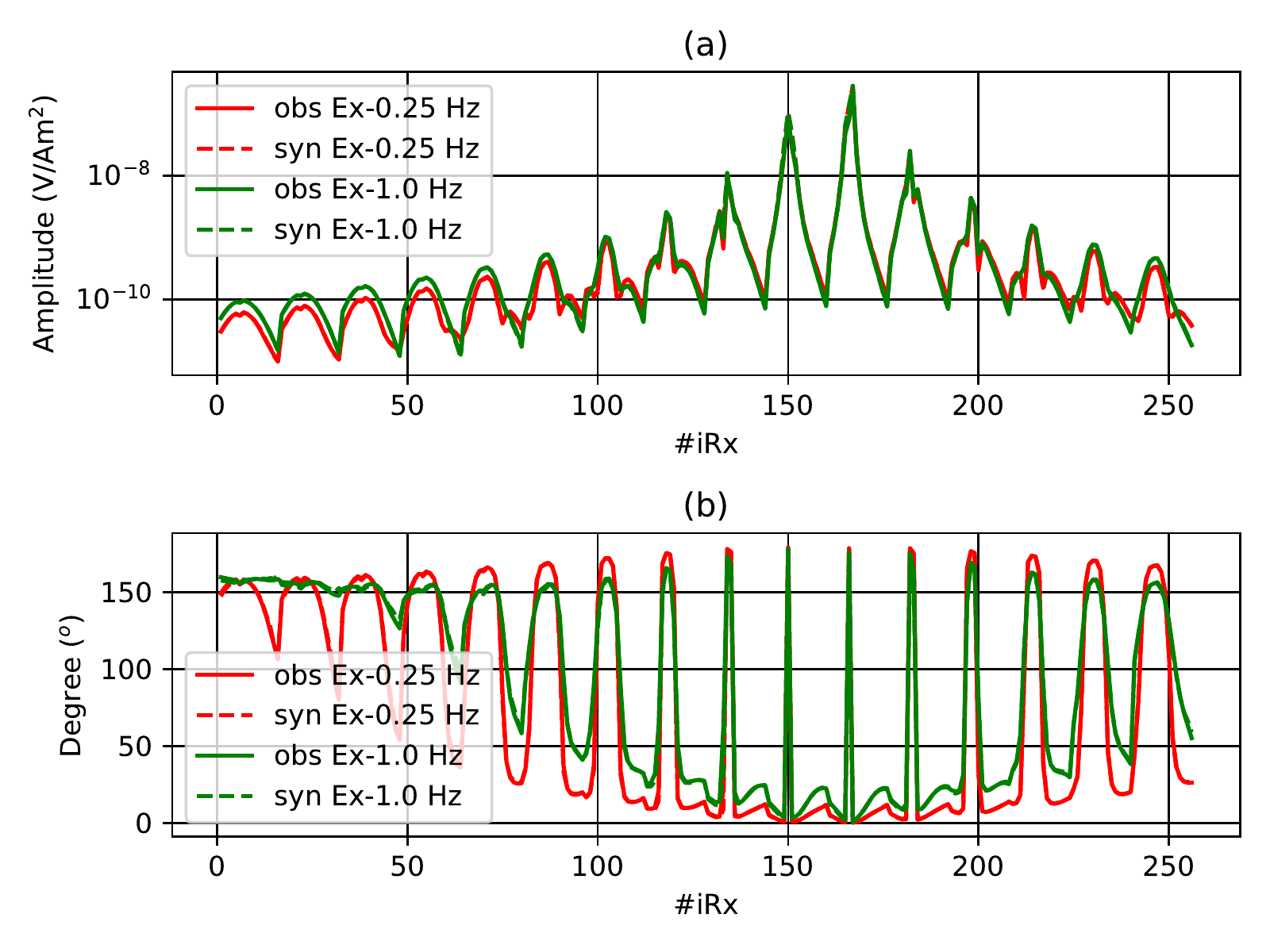}
  \caption{Comparison between observed data and synthetic data from inverted model according to (a) amplitude and (b) phase for transmitter Tx-10.  The synthetic data are now better aligned with the observed data after overlapping display.}\label{fig:dataiter30}
\end{figure}

\begin{figure}
  \centering
  \includegraphics[width=0.8\linewidth]{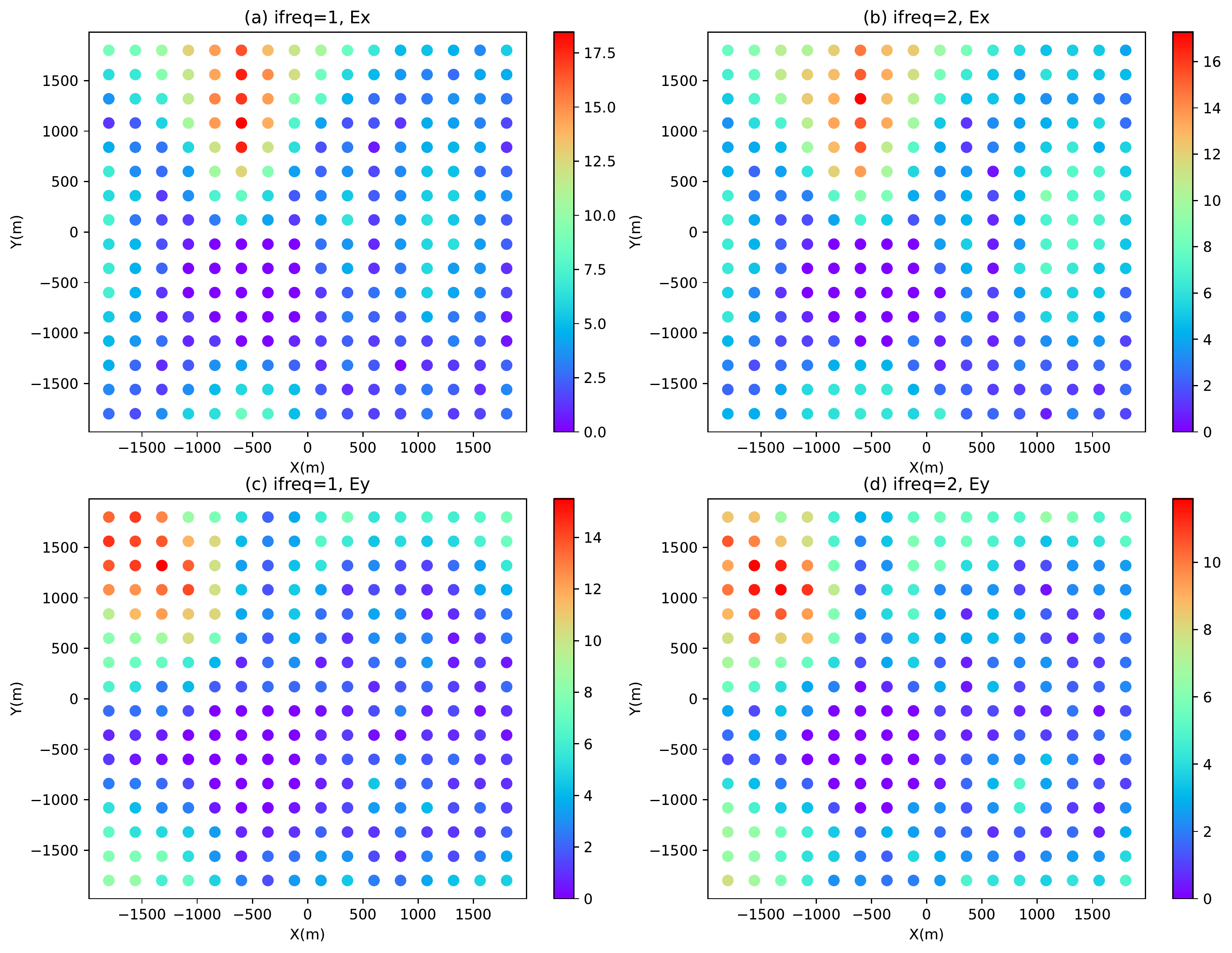}
  \caption{Scatter plot of the significant misfit for transmitter Tx-10 at iteration 1.}\label{fig:scatter_iter1}
\end{figure}

\begin{figure}
  \centering
  \includegraphics[width=0.8\linewidth]{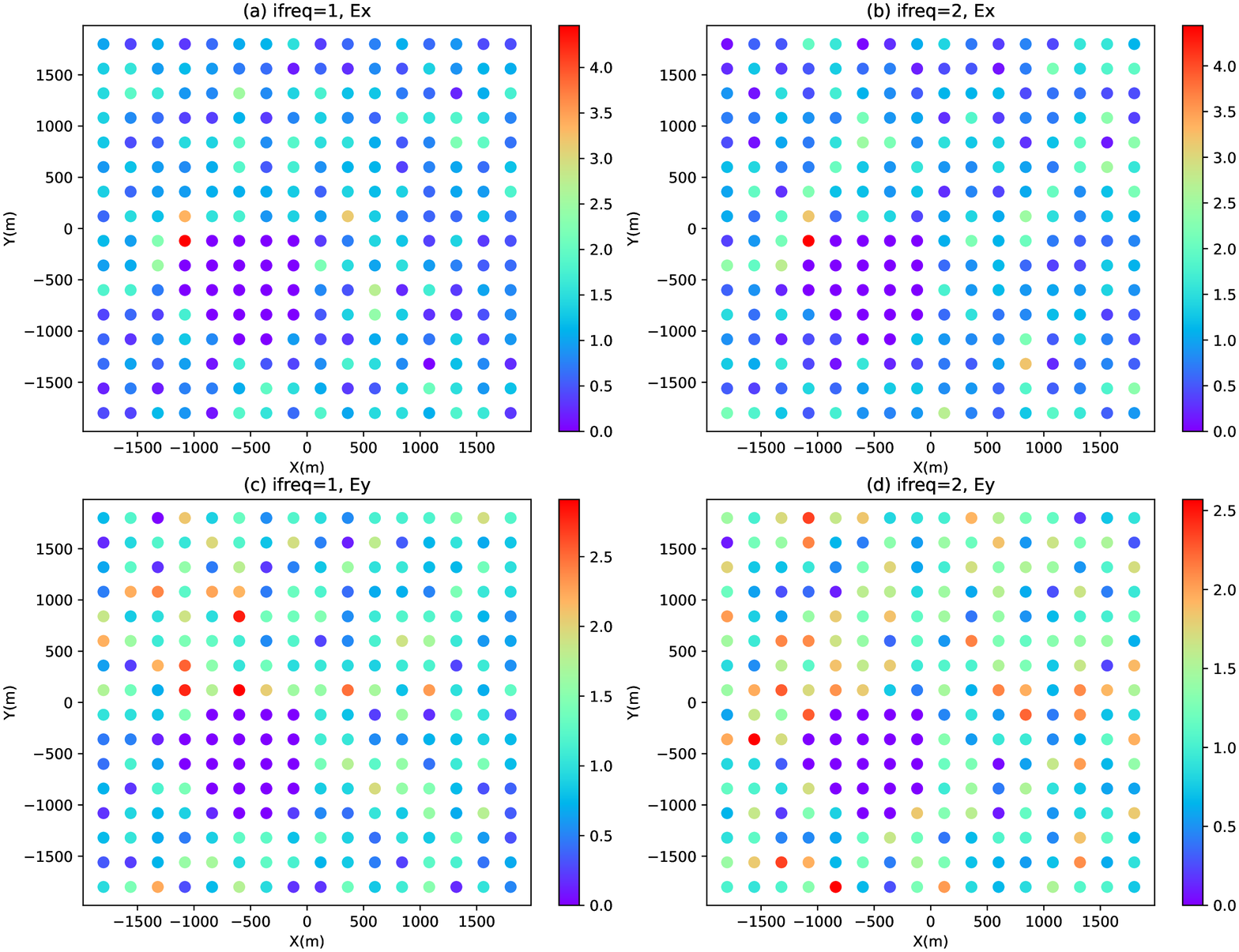}
  \caption{Scatter plot of the significant misfit for transmitter Tx-10 at iteration 30. The significant misfit becomes much lower after inversion.}\label{fig:scatter_iter30}
\end{figure}

This 3D CSEM inversion runs for 30 iterations. In order to check how the synthetic data fits the observations, we plot both the amplitude and the phase of the synthetic data at all receivers together with the observed data using the initial model and the retrieved resistivity model after inversion. In Figure~\ref{fig:dataiter1}, the synthetic data created from the homogeneous initial model cannot match the observed data very well. In Figure~\ref{fig:dataiter30}, the synthetic data generated from inverted resistivity volume indeed are better aligned with the observed data. Note that the horizontal axis in Figures~\ref{fig:dataiter1} and \ref{fig:dataiter30} are receiver indices rather than offset. This is due to  the broadside survey configuration. To better visualize the data matching connected to acquisition geometry, we design a scatter plot with hot colors indicating the data error at each receiver location. Figures~\ref{fig:scatter_iter1} shows that at the beginning of the inversion, the significant misfit (defined as $\|\mathbf{W}(\mathbf{d}_{obs}-\mathbf{d}_{syn})\|$ at each receiver location) is very large.  Since we mute the data within 600 m offset, the significant misfit surrounding transmitter Tx-10 are zeros. The plot in Figure~\ref{fig:scatter_iter1} also highlights the importance of azimuthal data has much larger significant misfit based than inline directions. After 30 inversion iterations, the significant misfit becomes much lower, see Figure~\ref{fig:scatter_iter30}.

\begin{figure}
  \centering
  \includegraphics[width=0.75\linewidth]{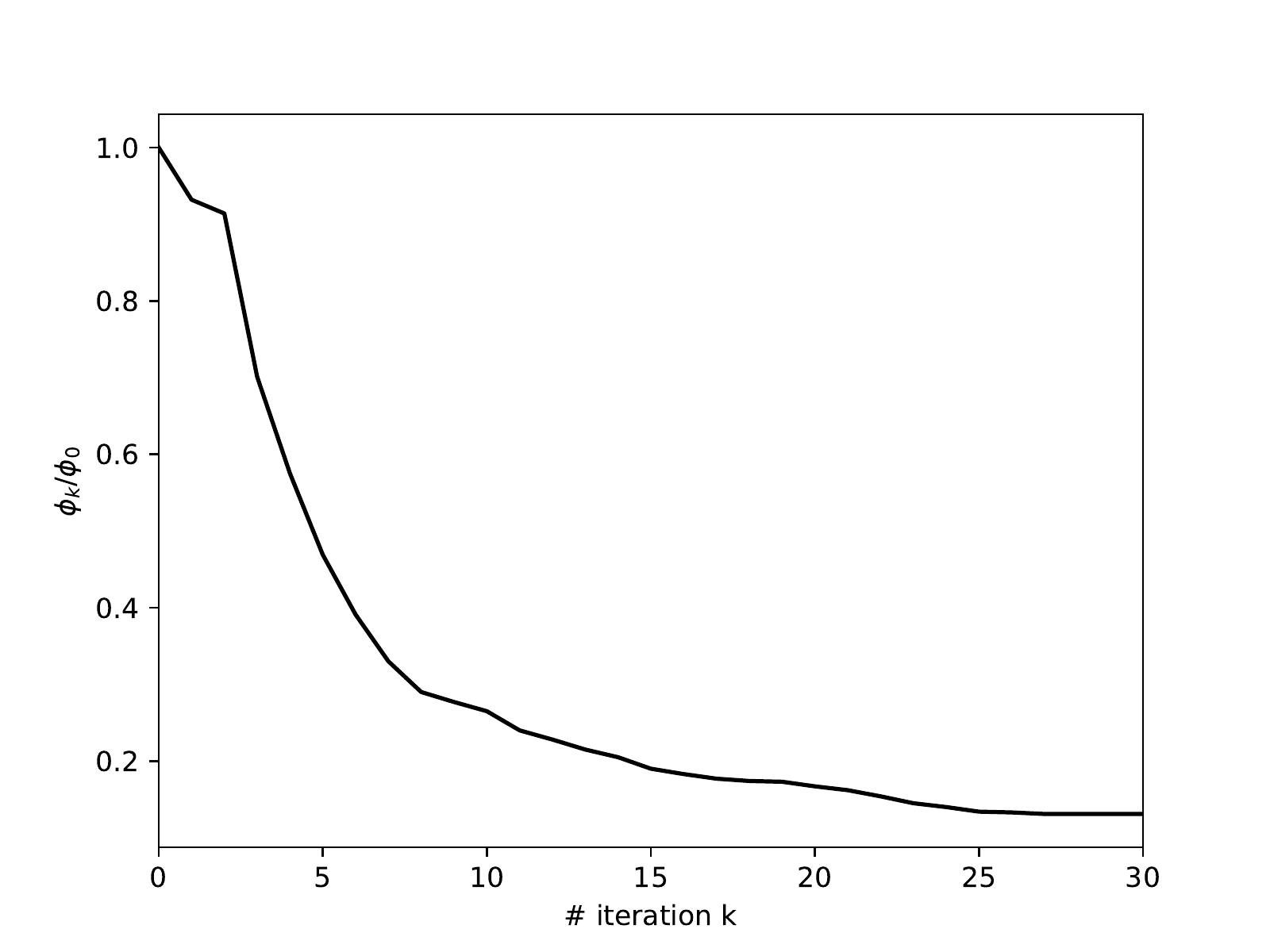}
  \caption{The convergence history in terms of the normalized misfit for the inversion of the two-block model.}\label{fig:convland}
\end{figure}

\begin{figure}
  \centering
  \includegraphics[width=0.7\linewidth]{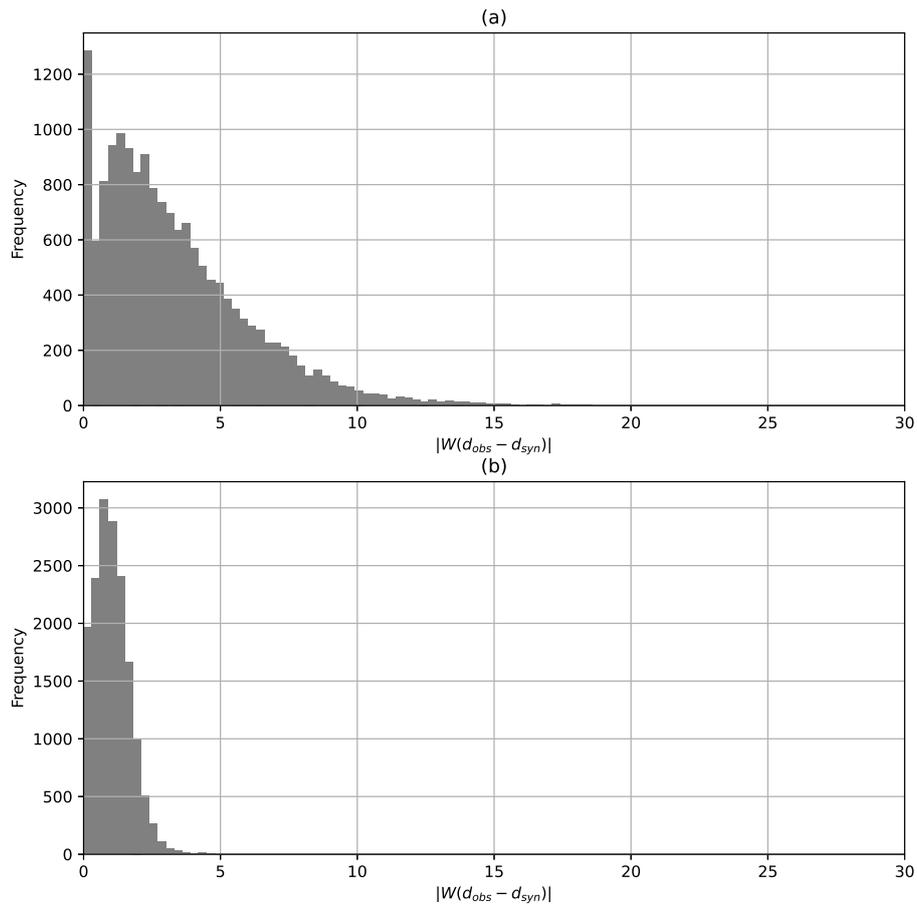}
  \caption{The histogram of the significant misfit  at (a) iteration 1 and (b) iteration 30.}\label{fig:landhistogram}
\end{figure}

\begin{figure}
  \centering
  \includegraphics[width=0.8\linewidth]{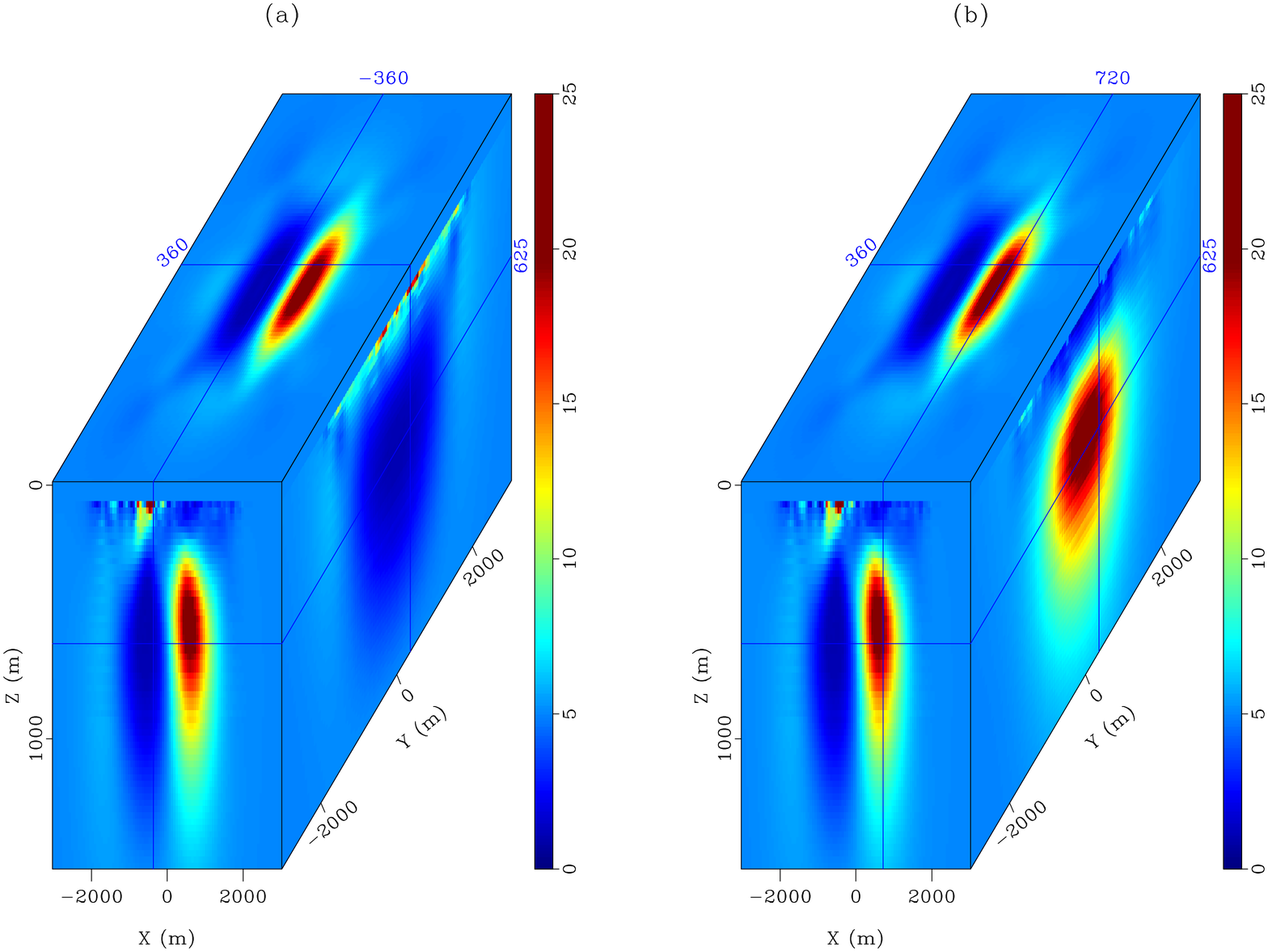}
  \caption{Inversion result for (a) low resistivity anomaly and (b) high resistivity anomaly.}\label{fig:landinv}
\end{figure}

The above data comparison gives us a confidence of our inversion scheme working properly, according to a specific source (Tx-10). In Figure~\ref{fig:convland}, we plot the normalized data misfit which is a global measure of the inversion. We see that after 30 iterations, it arrives at a relatively low normalized misfit. According to Figure~\ref{fig:landhistogram}a, the frequency of large significant misfit at receiver locations are distributed in a large range from 0 to 15. After the inversion, it has been compressed a lot (most of them are clustered within 3), as shown in Figure~\ref{fig:landhistogram}b.
The final inverted model is capable to retrieve both the low resistivity inclusion (Figure~\ref{fig:landinv}a) and the high resistivity inclusion (Figure~\ref{fig:landinv}b).

\subsection{A marine CSEM example}

\begin{figure}
  \centering
  \includegraphics[width=0.8\linewidth, angle=-90]{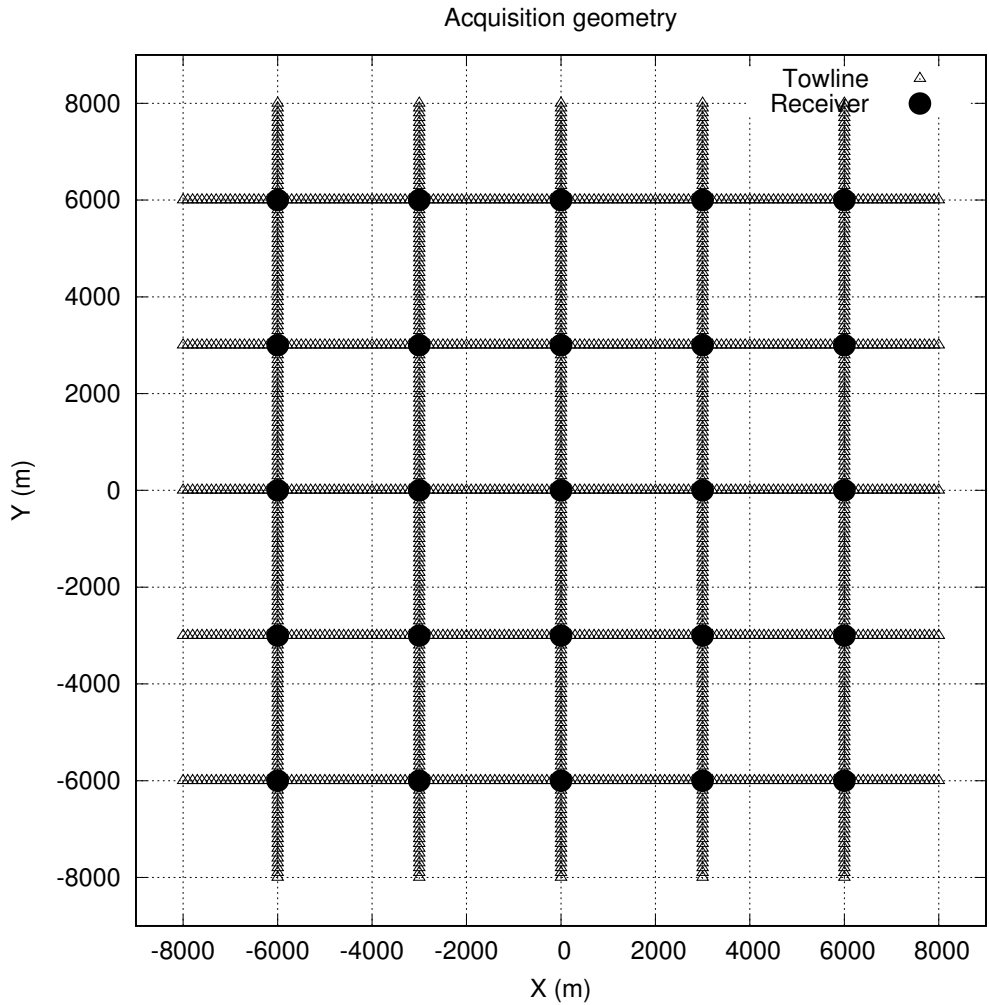}
  \caption{Survey layout sheet for marine CSEM inversion.}\label{fig:marinesurvey}
\end{figure}

In practical CSEM survey for reservoir exploration, a number of towlines using very low-frequency (0.1-10 Hz), high-energy electrical source will be deployed with hundreds of receivers. The source-receiver offset extends up to 10 km is very standard. Here,  a 3D CSEM inversion is carried out for 30 iterations based on a model of size $X=X_1\times X_2 \times X_3$, where $X_1=X_2=[-9000, 9000]$ m and $X_3=[0,3500]$ m. Figure~\ref{fig:marinesurvey} illustrates the survey layout sheet for the test: 10 towlines (5 in x direction and 5 in y direction) are deployed in order to cover the region of interest. Each towline has 81 source locations with a separation of 200 m. The receivers (25 in total) are deployed at the crossings of these towlines. We consider the reciprocity to switch the source and receiver, to achieve efficiency in inversion.

\begin{figure}
  \centering  
  \includegraphics[width=0.9\linewidth]{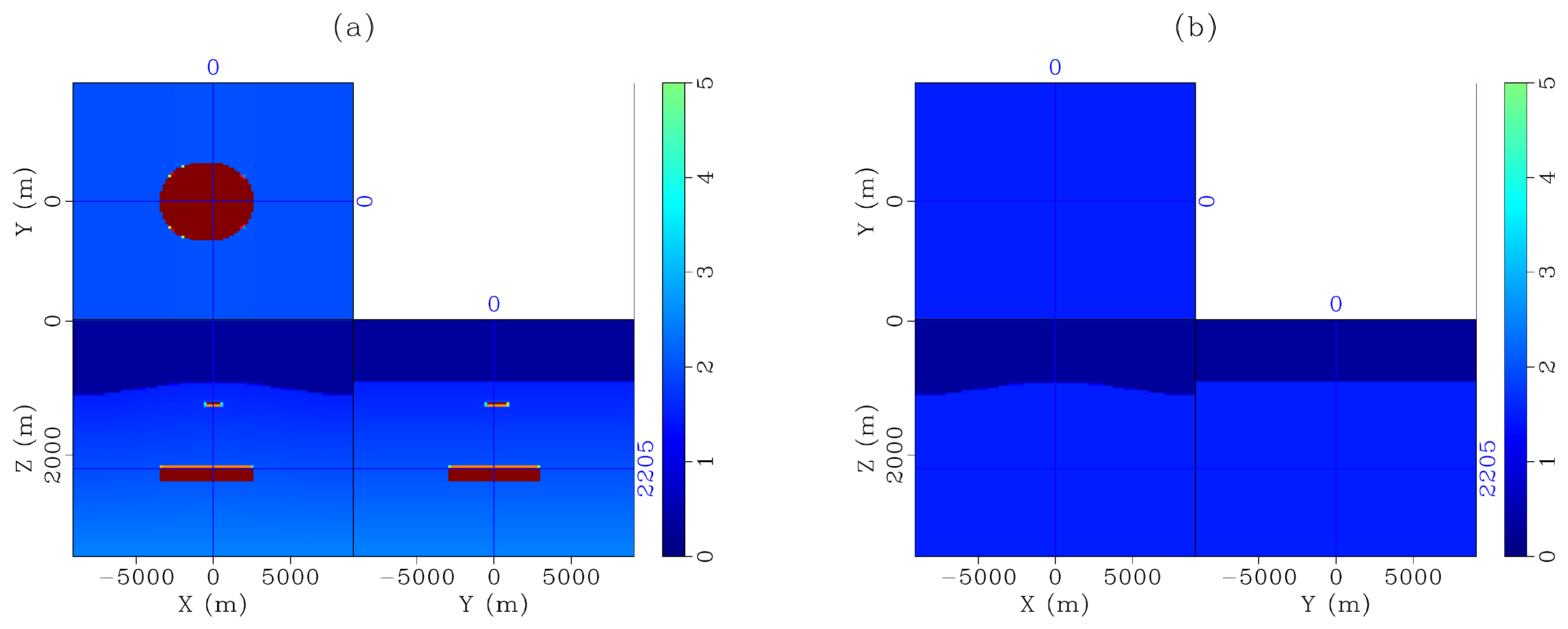}
  \caption{(a) The true resistivity model (note that the resistivity of the sediment below the seabed is varying from 1.5 to 2.5 $\Omega\cdot$m); (b) The initial resistivity model (the resistivity of the sediment takes a homogeneous value 1.5, to mimic the situation that no accurate a priori information is known). The models are clipped at 5 $\Omega\cdot$m for display purpose.}\label{fig:marinemodel}
\end{figure}

\begin{figure}
  \centering
  \includegraphics[width=0.8\linewidth]{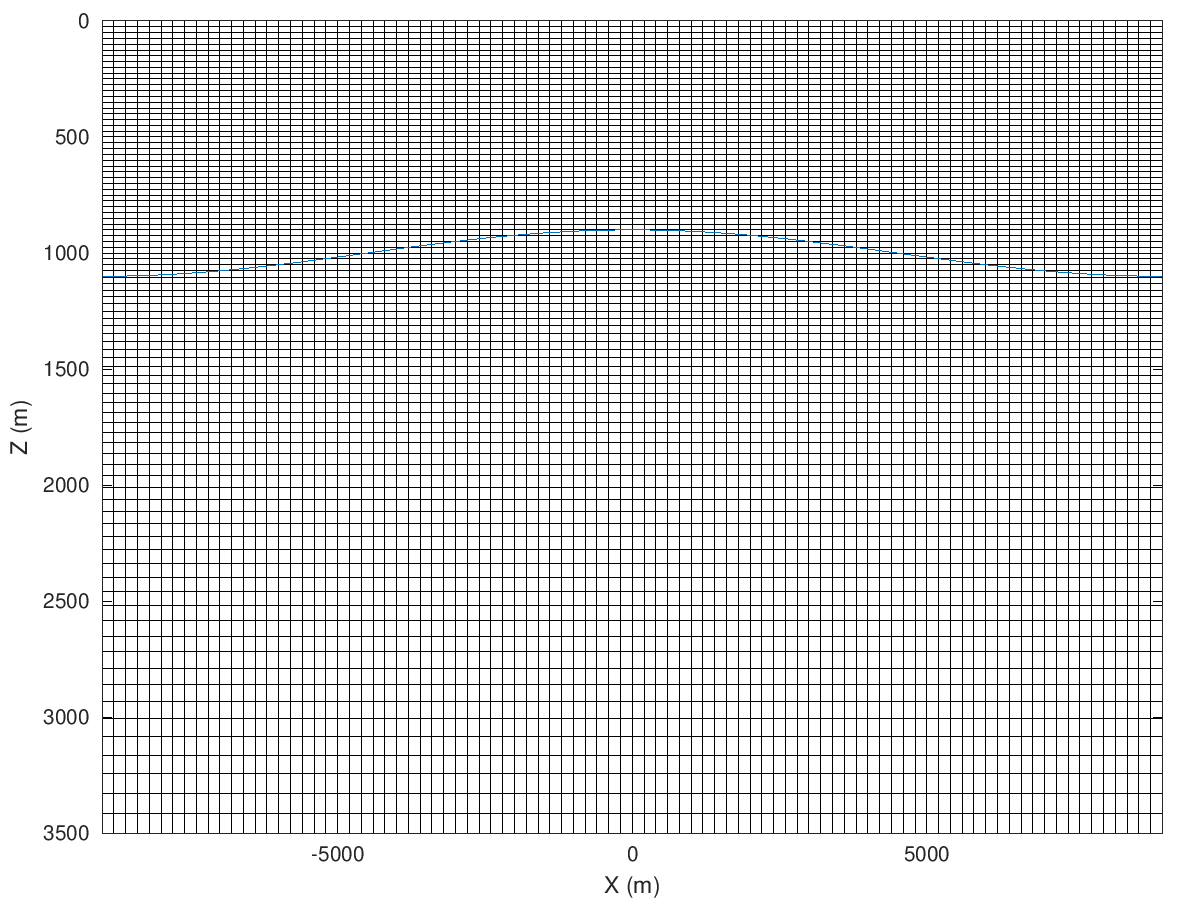}
  \caption{The vertical section of the nonuniform grid to mesh the 3D marine CSEM resistivity}\label{fig:mesh}
\end{figure}

As shown in Figure~\ref{fig:marinemodel}a, this synthetic model includes the seawater of 0.3 $\Omega\cdot$m. The resistivity of the background sediment is varying from 1.5 to 2.5 $\Omega\cdot$m along the depth. The most striking feature is a resistor of 10 $\Omega\cdot$m at shallow part and a disk-shaped resistor of 100 $\Omega\cdot$m sitting at the depth between 2200 m and 2350 m, which mimics the canonical reservoir. The starting model for the inversion in Figure~\ref{fig:marinemodel}b takes a homogeneous value 1.5 $\Omega\cdot$m to mimic the situation that no accurate a priori information is known for inversion.
A seafloor bathymetry with varying depth has been embedded in the model. The model is densely gridded around bathymetry and above, with the grid spacing $\Delta z=25$ m, see Figure~\ref{fig:mesh}. The grid has been stretched with a constant growing factor at a certain distance below the seabed.

\begin{figure}
  \centering
  \includegraphics[width=0.75\linewidth]{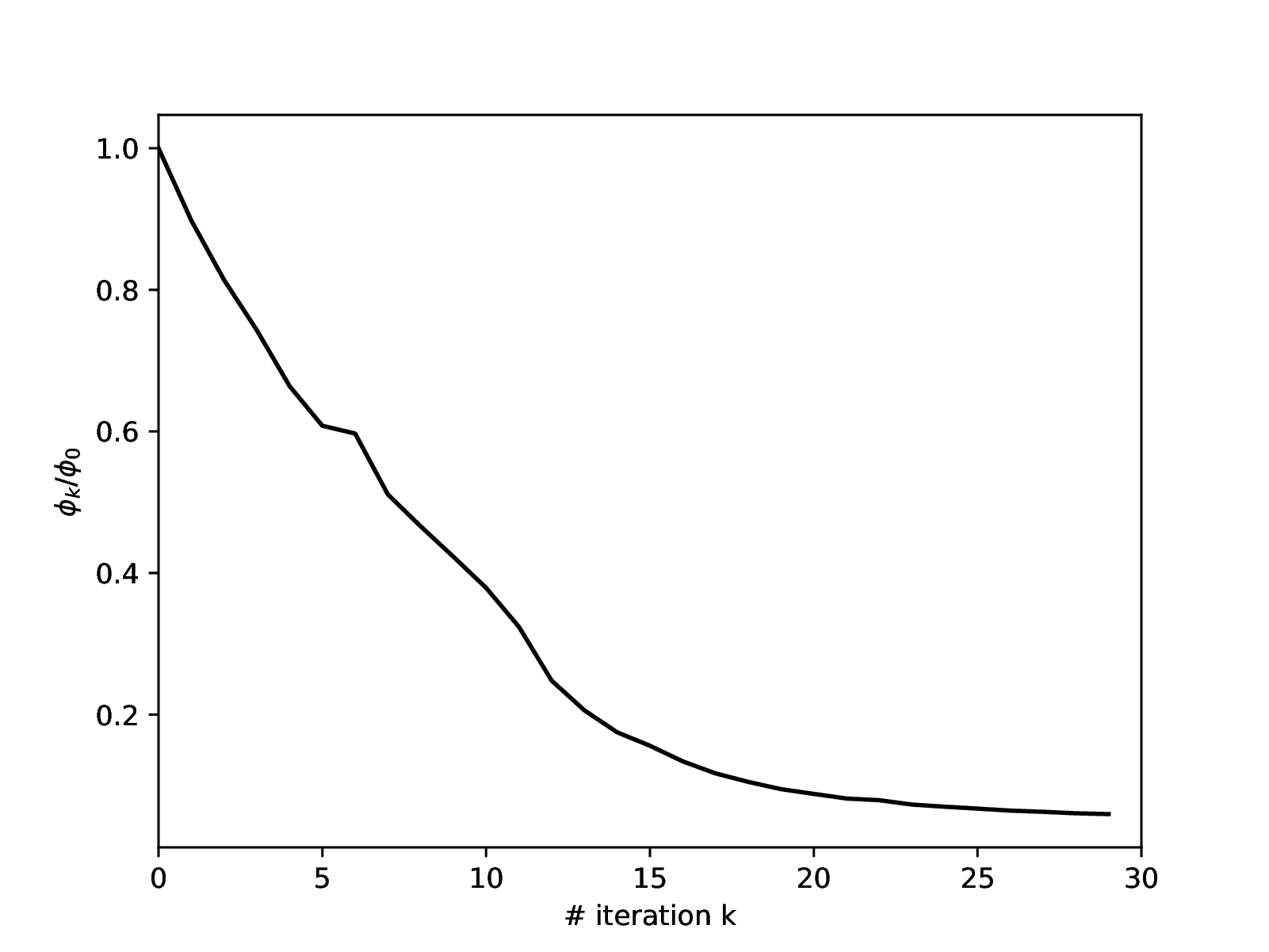}
  \caption{The convergence history of the normalized misfit for marine CSEM inversion.}\label{fig:convmarine}
\end{figure}

\begin{figure}
  \centering
  \includegraphics[width=\linewidth]{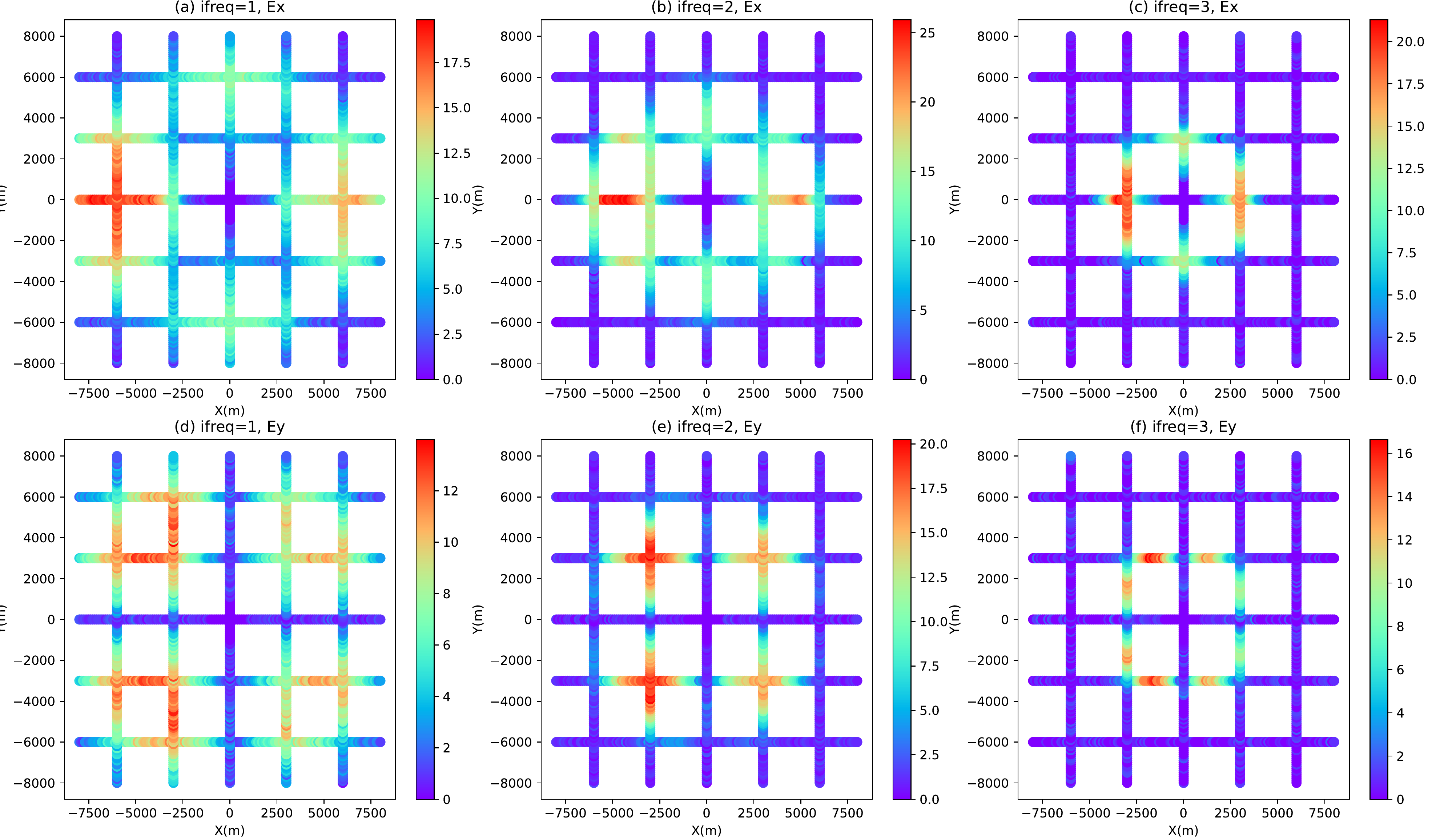}
  \caption{Scatter plot of the significant misfit for receiver at (0, 0) m in iteration 1.}\label{fig:mcsem_scatter_iter1}
\end{figure}

\begin{figure}
  \centering
  \includegraphics[width=\linewidth]{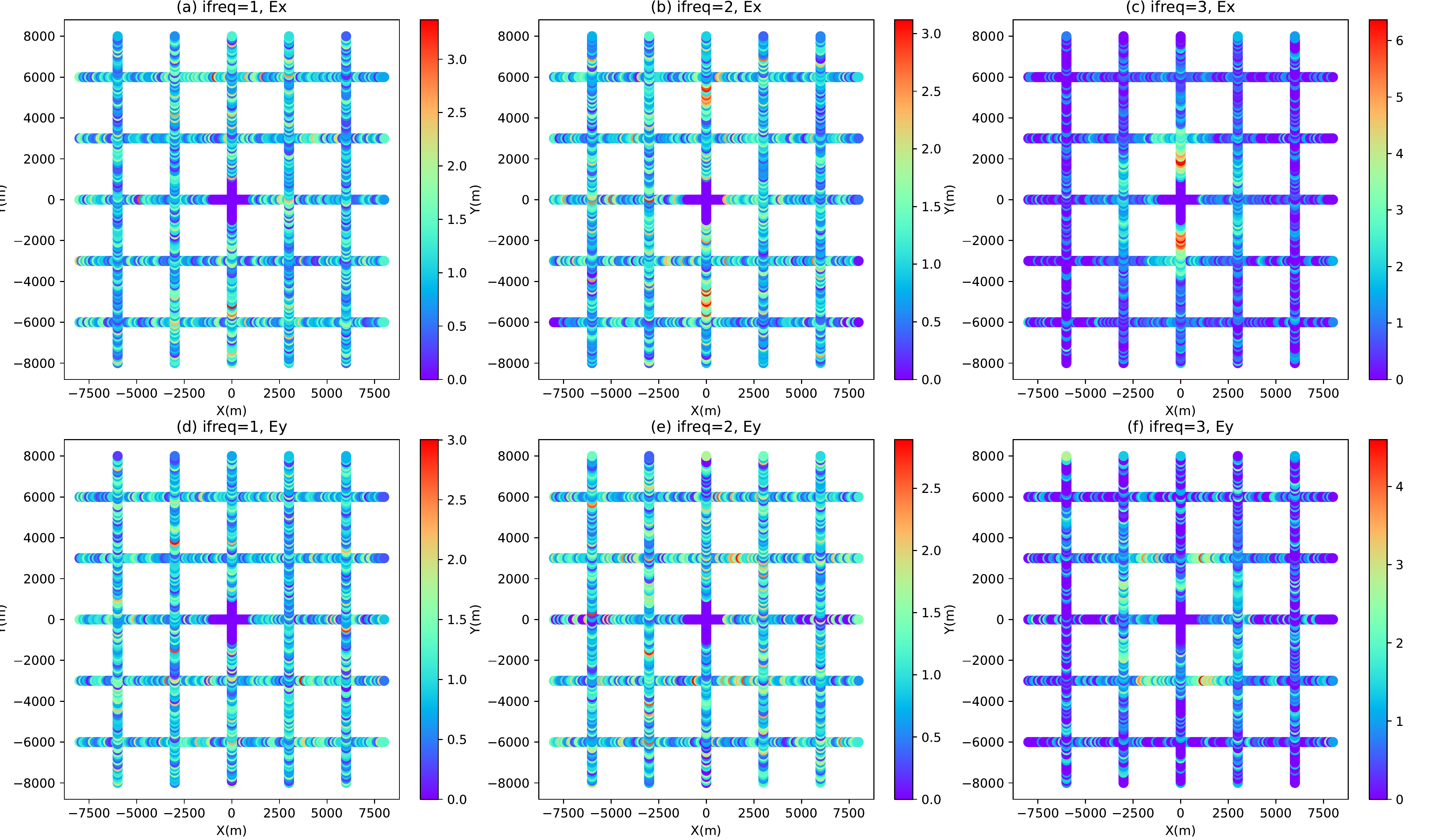}
  \caption{Scatter plot of the significant misfit for receiver at (0, 0) m in iteration 30. The significant misfit becomes much lower after inversion.}\label{fig:mcsem_scatter_iter30}
\end{figure}

\begin{figure}
  \centering
  \includegraphics[width=0.7\linewidth]{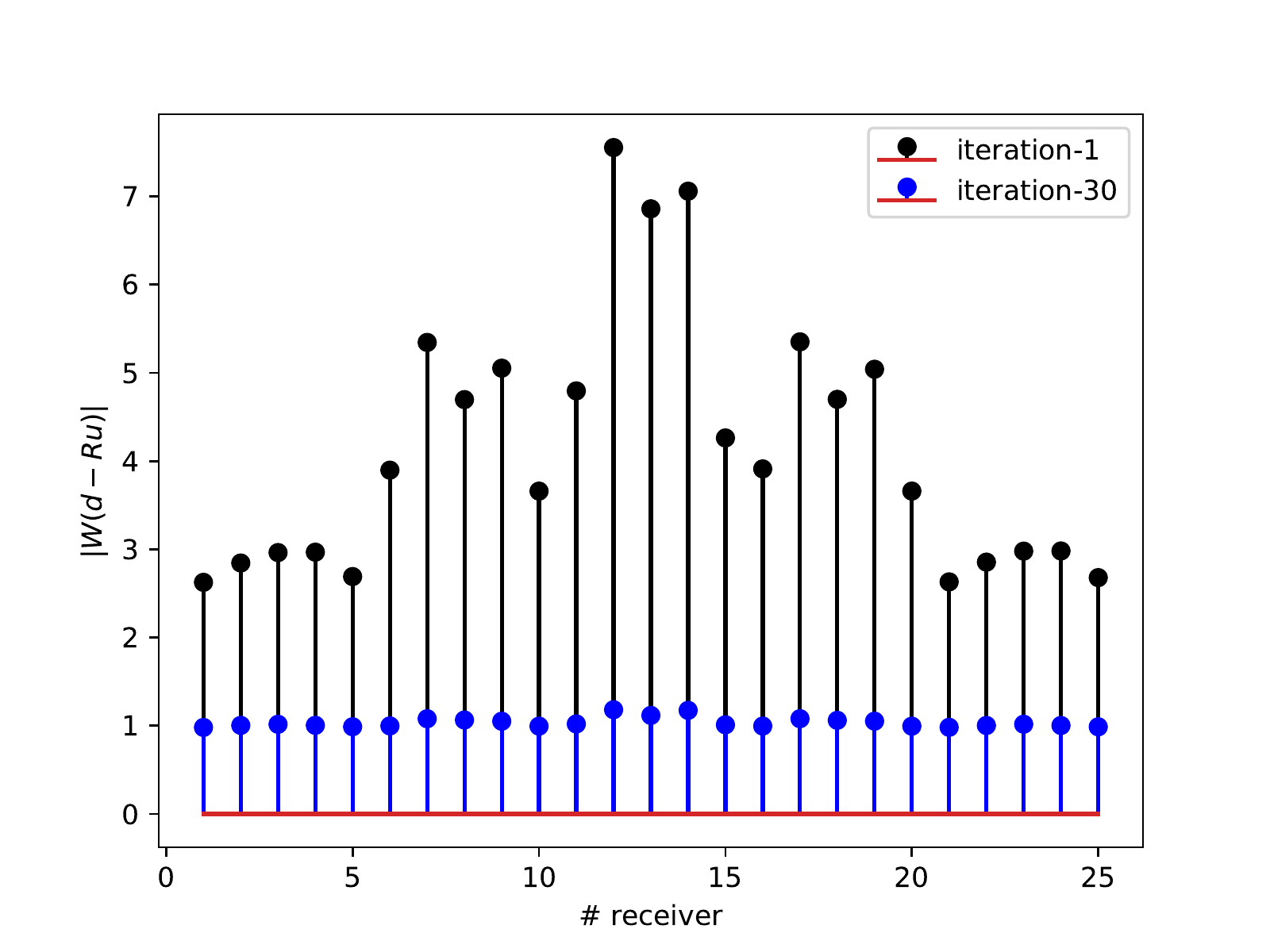}
  \caption{The significant misfit at each receiver location before and after inversion}\label{fig:rmse}
\end{figure}

\begin{figure}
  \centering
  \includegraphics[width=0.9\linewidth]{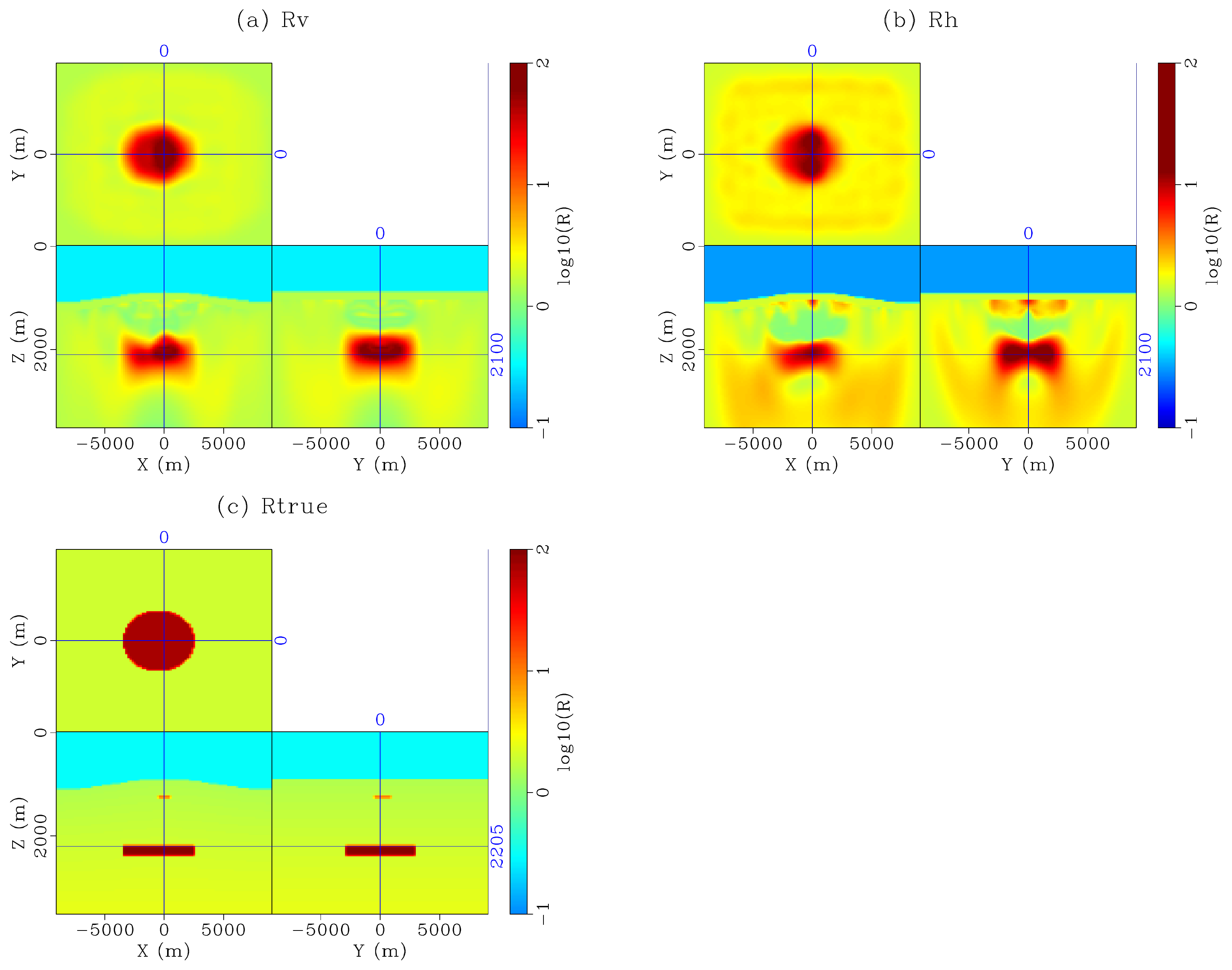}
  \caption{The result of 3D marine CSEM inversion: (a) vertical resistivity $\rho_v$, (b) horizontal resistivity $\rho_h$ and (c) true resistivity model. The display is in logorithmic scale. Note that in the inverted vertical resistivity, there is a small uphill variation just above the deep resistor which has a good correspondence to the shallow resistor in horizontal direction, where the highest resistivity has been retrieved at the center.}\label{fig:marine_inv}
\end{figure}

Three frequencies (0.25 Hz, 0.75 Hz and 2.25 Hz) have been used to perform this inversion. Figure~\ref{fig:convmarine} shows that this CSEM inversion was converging well.
The comparison between Figure~\ref{fig:mcsem_scatter_iter1} and Figure~\ref{fig:mcsem_scatter_iter30} shows that the significant misfit has been largely reduced for the receiver at the location (0, 0) m. Note that in Figure~\ref{fig:mcsem_scatter_iter30}c, the data misfit for the third frequency at the near offset is still large, indicating the presence of shallow resistor, which has not been well recovered in the inverted model (cf. Figure~\ref{fig:marine_inv}). To recover also the shallow resistor, we should use more near offset data rather than drop off them since 1000 m. A frequency dependent weighting strategy can be applied to boost the importance of high frequencies, as they have shallow penetrating depth according to skin depth. Since the scatter plot of the signficant misfit does not reflect the global misfit, the significant misfit at each receiver location has been shown in Figure~\ref{fig:rmse}: at the beginning of the inversion, they are far from 1; after the inversion, the significant misfit at all receiver locations are close to 1, which is the target value after adding noise into the observed data.

After 30 l-BFGS iterations, the inversion was successful to recover the resevoir in both vertical resistivity ($\rho_v=1/\sigma_h$) and the horizontal resistivity ($\rho_h=1/\sigma_v$), but the shallow resistor is difficult to obtain. The imprint of the shallow resistor can be found just above the resevoir in our reconstructed model: the horizontal location is correct but the depth is obviously incorrect.   It can be seen that the recovered anomaly in $\rho_v$ is of 100-150 m shallower than the true location.
This highlights the low resolution of CSEM inversion and the ambuiguity of the depth. 
The inverted vertical resistivity in Figure~\ref{fig:marine_inv}a is much better resolved than the horizontal resistivity in  Figure~\ref{fig:marine_inv}b. The value of the retrieved $\rho_v$ is higher than the retrieved $\rho_h$, while both are smaller than 100 $\Omega\cdot$m and thicker than the truth.  This is due to the effect of anomalous transverse resistance (ATR) \citep[equation 11]{mittet2012detection}: because the CSEM response is proportional to the product of resistivity difference times the thickness, lower resistivity with larger thickness can create similar amplitude response compared with higher resistivity with thin depth distribution. There is also weak increase of the resistivity in $\rho_h$ in  Figure~\ref{fig:marine_inv}b, suggesting that the initial model of homogeneous background lower than the true model, as can be seen by comparing with  Figure~\ref{fig:marine_inv}c.


\section{Discussion}

In our method, the modelling jobs for each source can be parallelized independently with perfect scalability. Since it is based on time stepping, the EM fields at the next time step will overwrite the memory unit at the previous step. This makes the memory usage of the method very economically. Our computation was performed on a Intel(R) Xeon(R) Gold 6258R CPU @ 2.70GHz, possessing 56 CPU cores (each core has 2 threads) and 128 GB memory. 
For the two-block model of 3D grid size 139*139*99 (after padding with PML boundaries on each side), combination of forward and adjoint modelling for building the gradient associated with one source only takes the memory of approximately 402 MB to compute all frequencies. (The memory requirement using frequency domain finite difference direct solver for such a large model size could not even fit 128 GB memory.)
This implies that we can simulateneously launch more than 100 modelling jobs (corresponding to 100 independent sources) to fully explore the small memory footprint of fictitous wave domain modelling in performing large scale 3D inversion. For 16 sources, we can use 4 threads to parallelize the computation by OpenMP using shared memory. The modelling takes less time at the beginning of the inversion, as the model is conductive and has less inhomogeneties. At the later iterations, each gradient building takes more computing time than before, because the resistive anomalies gradually come into the model, the time step $\Delta t$ becomes smaller and the number of timesteps becomes larger due to stability condition. The inversion with 30 iterations can be completed within 10 hours. On average, building one gradient (equivalent to one forward plus one adjoint solve) takes roughly 20 minutes.

It should be noted that the time series generated by our method only allows the matching of the field at given frequencies. The above showcases simply demonstrate that the proposed method works well, making fictitious wave domain an attractive modelling engine for 3D practical inversion. However, one has to be aware that the choice of regularization parameters $\beta$, $\alpha_x$, $\alpha_y$, $\alpha_z$ are of paramount importance to the final inversion result.   \citet{Lelievre2013gradient} have made a thorough investigation and concluded  that regularization of noisy data usually requires a lot of apriori knowledge about the geologic structure, such that the related parameters can be specified properly.  Our choices may be far from perfect, and better choices should be possible to significantly improve the inversion result, while avoiding the false positives in imaging.  These are, however, irrelevant to the core contribution of this work.

The derivation in the appendix \ref{sec:appendix} develops an alternative adjoint formulation to avoid solving the regularized linear inverse problem. Equation \eqref{eq:chainrule}
gives an explicit expression to build CSEM inversion gradient directly by cross-correlation of fictitious time domain forward and adjoint fields. However, 3D implementation of such a strategy is not trivial as it is mandatory to simultaneously access fictitious time domain forward and adjoint field at each time step. This is challenging due to the opposite directions in time between forward and adjoint fields. In doing so, one may think of storing all 3D wavefield snapshots during thunsands of times and accessing them one by one in reverse order, which is extremely inefficient due to slow memory acess speed and heavy IO traffic. A feasible approach is to use the wavefield reconstruction method by storing the boundary values \citep{Yang_2016_WRB}, or resort to optimal checkpointing strategy \citep{Symes_2007_RTM,Yang_2016_CAR}. Wavefield reconstruction by boundary storage implies each gradient computation involves three times of simulation, two for forward fields and one for adjoint fields. Optimal checkpointing requires much higher recomputation ratio and thus even more computationally expensive. In summary, these strategies demand either significant memory consumption for storing boundaries of forward fields at every time step, or  sophisticated operations to achieve optimal combinatorial recomputation of  the forward fields, thus left as the future work. The scheme used in this paper is a hybrid combination of time domain modelling and frequency domain inversion. To build the gradient, only two times of modelling are needed, one for forward fields and one for adjoint fields. Consequently, the gradient computed by product between forward and adjoint fields in frequency domain becomes the most efficient option among all equivalent formulations.

\section{Conclusion}

A nonlinear CSEM inversion based on fictitious wave domain modelling using timestepping simulation has been presented. To construct the gradient using forward and adjoint fields, the adjoint source time functions are estimated by solving a regularized linear inverse problem. A matrix-free implementation using the basis functions is proposed such that the linear inverse problem is iteratively solved only  once and then the resulting basis functions can be reused to build up all adjoint source time functions at every receiver location.
The method enables efficient 3D CSEM inversion which are demonstrated by numerical examples.
The method will become more attractive for land CSEM which normally requires much more number of frequencies than marine CSEM. The proposed method can equally be applied to magnetotelluric (MT) inversion to decipher the subsurface structure in even larger scale using much lower frequencies of the EM recordings.

\section*{Computer Code Availability}

Name of the code: \verb|libEMMI|

Functionality: Library for 3D controlled-source electromagnetic modelling and inversion

Contact: Pengliang Yang, Harbin Institute of Technology, China

E-mail: ypl.2100@gmail.com

System requirements: Linux OS

Programming language: C, Fortran

Result visualization: Python3, Madagascar (Scons) and Gnuplot

Compilation requirement: gcc compiler, mpicc compiler and make

Software required: fftw3 (\url{http://fftw.org/})

The code has been released publically via github repository \url{https://github.com/yangpl/libEMMI}.

\section*{Acknowledgements}

The author acknowledges the support from Chinese Fundamental Research Funds for the Central
Universities (AUGA5710010121) and National Natural Science Fundation of China (42274156).
The author is indebted to Rune Mittet for the inspiring discussions and the encouragement of this work. The alternative adjoint formulation in the appendix \ref{sec:appendix} was due to Ren\'e-\'Edourt Plessix.

\appendix

\section{An alternative adjoint formulation}\label{sec:appendix}

The above approach uses diffusive Maxwell equation in frequency domain as the state equation, which holds no matter what kind of modelling engine is used as the solver. 
The misfit in \eqref{eq:phid} can also be expanded as
\begin{displaymath}
  \phi_d(m)=\frac{1}{2}\sum_s\sum_k|\mathbf{W}_1(E^{obs}(\mathbf{x}_r,\omega_k;\mathbf{x}_s)-RE(\mathbf{x},\omega_k;\mathbf{x}_s))|^2 +
  |\mathbf{W}_2(H^{obs}(\mathbf{x}_r,\omega_k;\mathbf{x}_s)-RH(\mathbf{x},\omega_k;\mathbf{x}_s))|^2.
\end{displaymath}
Based on equation \eqref{eq:fictitiouswave} and the correspondence between wave and diffusive EM fields, the Lagrangian functional translates a constrained minimization into an unconstrained optimization
\begin{equation}
  \begin{split}
    \mathcal{L} = &\frac{1}{2}\sum_s\sum_k|\mathbf{W}_1(E^{obs}(\mathbf{x}_r,\omega_k;\mathbf{x}_s)-R E(\mathbf{x},\omega_k;\mathbf{x}_s))|^2 +
  |\mathbf{W}_2(H^{obs}(\mathbf{x}_r,\omega_k;\mathbf{x}_s)-R H(\mathbf{x},\omega_k;\mathbf{x}_s))|^2 \\
    +&\Re\sum_s\sum_k \int_X \mathrm{d}\mathbf{x} \overline{\lambda_1(\mathbf{x},\omega_k;\mathbf{x}_s)} \bigg(E(\mathbf{x},\omega_k;\mathbf{x}_s)-\int_0^{T_{\max}} dt' E'(\mathbf{x},t';\mathbf{x}_s) e^{-\sqrt{\omega_k\omega_0}t'} e^{\mathrm{i}\sqrt{\omega_k\omega_0}t'}\bigg)\\
    +&\Re\sum_s\sum_k \int_X \mathrm{d}\mathbf{x} \overline{\lambda_2(\mathbf{x},\omega_k;\mathbf{x}_s)}\bigg( H(\mathbf{x},\omega_k;\mathbf{x}_s)-\sqrt{-\frac{2\omega_0}{\mathrm{i}\omega_k}}\int_0^{T_{\max}} dt' H'(\mathbf{x},t';\mathbf{x}_s) e^{-\sqrt{\omega_k\omega_0}t'} e^{\mathrm{i}\sqrt{\omega_k\omega_0}t'}\bigg)\\
    +&\sum_s\int_0^{T_{\max}} dt'\int_X \mathrm{d}\mathbf{x} \underline{H}'
    (\mathbf{x},t';\mathbf{x}_s)\cdot\bigg(\nabla\times E'(\mathbf{x},t';\mathbf{x}_s) + \mu'\partial_t H'(\mathbf{x},t';\mathbf{x}_s) -M'(\mathbf{x}_s,t')\bigg)\\
    +&\sum_s\int_0^{T_{\max}} dt'\int_X \mathrm{d}\mathbf{x} \underline{E}'(\mathbf{x},t';\mathbf{x}_s)\cdot\bigg(\nabla\times H'(\mathbf{x},t';\mathbf{x}_s) - \epsilon\partial_t E'(\mathbf{x},t';\mathbf{x}_s) -J'(\mathbf{x}_s,t')\bigg),
  \end{split}
\end{equation}
where $\lambda_1(\mathbf{x},\omega;\mathbf{x}_s)$ and $\lambda_2(\mathbf{x},\omega;\mathbf{x}_s)$ are Lagrangian multipliers in the frequency domain, while  $\underline{E}'(\mathbf{x},t';\mathbf{x}_s)$ and $\underline{H}'(\mathbf{x},t';\mathbf{x}_s)$ are Lagrangian multipliers in the fictitious time domain.

Setting the Lagrangian with respect to state variables to zero ($\partial\mathcal{L}/\partial \lambda_1=\partial\mathcal{L}/\partial \lambda_2=\partial\mathcal{L}/\partial \underline{E}'=\partial \mathcal{L}/\underline{H}'=0$) gives exactly the state equations. The adjoint equations are obtained by setting the Lagrangian with respect to state variables, i.e.,
\begin{subequations}
  \begin{align}
  \frac{\partial\mathcal{L}}{\partial E(\mathbf{x},\omega_k;\mathbf{x}_s)}=0 &
  \Leftrightarrow
  \lambda_1(\mathbf{x},\omega_k;\mathbf{x}_s) = R^\dagger \mathbf{W}_1^\dagger \mathbf{W}_1(E^{obs}(\mathbf{x}_r, \omega_k)-R E(\mathbf{x},\omega_k;\mathbf{x}_s)),\label{eq:lambda1}\\
  \frac{\partial\mathcal{L}}{\partial H(\mathbf{x},\omega_k;\mathbf{x}_s)}=0 &
  \Leftrightarrow
  \lambda_2(\mathbf{x},\omega_k;\mathbf{x}_s) = R^\dagger \mathbf{W}_2^\dagger \mathbf{W}_2(H^{obs}(\mathbf{x}_r, \omega_k;\mathbf{x}_s)-R H(\mathbf{x},\omega_k;\mathbf{x}_s)),\label{eq:lambda2}\\
  \frac{\partial\mathcal{L}}{\partial E'(\mathbf{x},t';\mathbf{x}_s)}=0 &
  \Leftrightarrow
  \nabla\times \underline{H}'(\mathbf{x},t';\mathbf{x}_s) + \epsilon\partial_t \underline{E}'(\mathbf{x},t';\mathbf{x}_s) =\Re \sum_k \overline{\lambda_1(\mathbf{x},\omega_k;\mathbf{x}_s)} e^{-\sqrt{\omega_k\omega_0}t'} e^{\mathrm{i}\sqrt{\omega_k\omega_0}t'},\label{eq:appadj1}\\
  \frac{\partial\mathcal{L}}{\partial H'(\mathbf{x},t';\mathbf{x}_s)}=0 &
  \Leftrightarrow
  -\mu\partial_t \underline{H}'(\mathbf{x},t';\mathbf{x}_s) + \nabla\times \underline{E}'(\mathbf{x},t';\mathbf{x}_s) = \Re\sum_k \sqrt{-\frac{2\omega_0}{\mathrm{i}\omega_k}}\overline{\lambda_2(\mathbf{x},\omega_k;\mathbf{x}_s)} e^{-\sqrt{\omega_k\omega_0}t'} e^{\mathrm{i}\sqrt{\omega_k\omega_0}t'}.\label{eq:appadj2}
  \end{align}
\end{subequations}
Equations \eqref{eq:lambda1} and \eqref{eq:lambda2} reveals that $\lambda_1$ and $\lambda_2$ are nothing more than the electric and magnetic components of the weighted data residual, while equations \eqref{eq:appadj1} and \eqref{eq:appadj2} form a new set of Maxwell system. It is important to remark that in deriving the above equations, we have tacitly applied integration by parts in time and space, assuming zero initial condition of the forward EM fields
\begin{equation}
  E'(\mathbf{x},t'=0)=H'(\mathbf{x},t'=0)= 0,\quad \mathbf{x} \in X,
\end{equation}
zero final conditions of the adjoint fields
\begin{equation}
  \underline{E}'(\mathbf{x},t'=T_{\max})=\underline{H}'(\mathbf{x},t'=T_{\max})=0,
  \quad \mathbf{x}\in X,
\end{equation}
and homogeneous boundary condition in space for both forward and adjoint fields
\begin{equation}
  E'(\mathbf{x},\cdot)=H'(\mathbf{x},\cdot)
  =\underline{E}'(\mathbf{x},\cdot)=\underline{H}'(\mathbf{x},\cdot)=0,
  \quad \mathbf{x} \in \partial X.
\end{equation}

At the saddle point, we again obtain 
\begin{equation}
  \frac{\partial \phi_d}{\partial \epsilon}=\frac{\partial \mathcal{L}}{\partial \epsilon}=-\sum_s\int_0^{T_{\max}} \mathrm{d}t' \underline{E}'(\mathbf{x},t';\mathbf{x}_s)\cdot\partial_t E'(\mathbf{x},t';\mathbf{x}_s).
\end{equation}
Since $\sigma_{ij}=2\omega_0\epsilon_{ij}$, the application of the chain rule yields
\begin{equation}\label{eq:chainrule}
\frac{\partial \phi_d}{\partial \sigma_{ij}}=\frac{\partial \phi_d}{\partial \epsilon_{ij}}\frac{\partial\epsilon_{ij}}{\partial \sigma_{ij}}=-\sum_s\frac{1}{2\omega_0}\int_0^{T_{\max}} \mathrm{d}t' \underline{E}'_j(\mathbf{x},t';\mathbf{x}_s)\cdot\partial_t E'_i(\mathbf{x},t';\mathbf{x}_s).
\end{equation}
With the gradient at hand, one can then construct the descent direction based on nonlinear minimization methods to solve the inverse problem.

\bibliographystyle{cas-model2-names}
\newcommand{\SortNoop}[1]{}


\begin{thebibliography}{55}
\expandafter\ifx\csname natexlab\endcsname\relax\def\natexlab#1{#1}\fi
\providecommand{\url}[1]{\texttt{#1}}
\providecommand{\href}[2]{#2}
\providecommand{\path}[1]{#1}
\providecommand{\DOIprefix}{doi:}
\providecommand{\ArXivprefix}{arXiv:}
\providecommand{\URLprefix}{URL: }
\providecommand{\Pubmedprefix}{pmid:}
\providecommand{\doi}[1]{\href{http://dx.doi.org/#1}{\path{#1}}}
\providecommand{\Pubmed}[1]{\href{pmid:#1}{\path{#1}}}
\providecommand{\bibinfo}[2]{#2}
\ifx\xfnm\relax \def\xfnm[#1]{\unskip,\space#1}\fi
\bibitem[{Abubakar et~al.(2008)Abubakar, Habashy, Druskin, Knizhnerman and
  Alumbaugh}]{abubakar20082}
\bibinfo{author}{Abubakar, A.}, \bibinfo{author}{Habashy, T.},
  \bibinfo{author}{Druskin, V.}, \bibinfo{author}{Knizhnerman, L.},
  \bibinfo{author}{Alumbaugh, D.}, \bibinfo{year}{2008}.
\newblock \bibinfo{title}{2.5 {D} forward and inverse modeling for interpreting
  low-frequency electromagnetic measurements}.
\newblock \bibinfo{journal}{Geophysics} \bibinfo{volume}{73},
  \bibinfo{pages}{F165--F177}.
\bibitem[{Alumbaugh and Newman(1997)}]{alumbaugh1997three}
\bibinfo{author}{Alumbaugh, D.}, \bibinfo{author}{Newman, G.},
  \bibinfo{year}{1997}.
\newblock \bibinfo{title}{{Three-dimensional massively parallel electromagnetic
  inversion-II. Analysis of a crosswell electromagnetic experiment}}.
\newblock \bibinfo{journal}{Geophysical Journal International}
  \bibinfo{volume}{128}, \bibinfo{pages}{355--363}.
\bibitem[{Bj{\"{o}}rck(1996)}]{Bjorck_1996_NML}
\bibinfo{author}{Bj{\"{o}}rck, {\AA{}}.}, \bibinfo{year}{1996}.
\newblock \bibinfo{title}{Numerical methods for least squares problems}.
\newblock \bibinfo{publisher}{SIAM, Society for Industrial and Applied
  Mathematics, Philadelphia}.
\bibitem[{Chang-Chun et~al.(2015)Chang-Chun, Xiu-Yan, Yun-He, Yan-Fu, Chang-Kai
  and Jing}]{yin2015review}
\bibinfo{author}{Chang-Chun, Y.}, \bibinfo{author}{Xiu-Yan, R.},
  \bibinfo{author}{Yun-He, L.}, \bibinfo{author}{Yan-Fu, Q.},
  \bibinfo{author}{Chang-Kai, Q.}, \bibinfo{author}{Jing, C.},
  \bibinfo{year}{2015}.
\newblock \bibinfo{title}{Review on airborne electromagnetic inverse theory and
  applications}.
\newblock \bibinfo{journal}{Geophysics} \bibinfo{volume}{80},
  \bibinfo{pages}{W17--W31}.
\bibitem[{Chave and Cox(1982)}]{chave1982on}
\bibinfo{author}{Chave, A.D.}, \bibinfo{author}{Cox, C.S.},
  \bibinfo{year}{1982}.
\newblock \bibinfo{title}{{Controlled Electromagnetic Sources for Measuring
  Electrical Conductivity Beneath the Oceans 1. Forward Problem and Model
  Study}}.
\newblock \bibinfo{journal}{Journal of Geophysical Research}
  \bibinfo{volume}{87}, \bibinfo{pages}{5327–5338}.
\bibitem[{Claerbout and Fomel(2008)}]{claerbout2008image}
\bibinfo{author}{Claerbout, J.F.}, \bibinfo{author}{Fomel, S.},
  \bibinfo{year}{2008}.
\newblock \bibinfo{title}{{Image estimation by example: geophysical soundings
  image construction: multidimensional autoregression}}.
\newblock \bibinfo{publisher}{Stanford University}.
\bibitem[{Commer and Newman(2008)}]{commer2008new}
\bibinfo{author}{Commer, M.}, \bibinfo{author}{Newman, G.A.},
  \bibinfo{year}{2008}.
\newblock \bibinfo{title}{New advances in three-dimensional controlled-source
  electromagnetic inversion}.
\newblock \bibinfo{journal}{Geophysical Journal International}
  \bibinfo{volume}{172}, \bibinfo{pages}{513--535}.
\bibitem[{Constable(2010)}]{constable2010}
\bibinfo{author}{Constable, S.}, \bibinfo{year}{2010}.
\newblock \bibinfo{title}{{Ten years of marine CSEM for hydrocarbon
  exploration}}.
\newblock \bibinfo{journal}{Geophysics} \bibinfo{volume}{75},
  \bibinfo{pages}{75A67--75A81}.
\bibitem[{Constable and Srnka(2007)}]{constable2007introduction}
\bibinfo{author}{Constable, S.}, \bibinfo{author}{Srnka, L.J.},
  \bibinfo{year}{2007}.
\newblock \bibinfo{title}{An introduction to marine controlled-source
  electromagnetic methods for hydrocarbon exploration}.
\newblock \bibinfo{journal}{Geophysics} \bibinfo{volume}{72},
  \bibinfo{pages}{WA3--WA12}.
\bibitem[{Constable et~al.(1986)Constable, Cox and
  Chave}]{constable1986offshore}
\bibinfo{author}{Constable, S.C.}, \bibinfo{author}{Cox, C.S.},
  \bibinfo{author}{Chave, A.D.}, \bibinfo{year}{1986}.
\newblock \bibinfo{title}{Offshore electro-magnetic surveying techniques}, in:
  \bibinfo{booktitle}{56th Annual International Meeting, SEG, Expanded
  Abstracts}, \bibinfo{publisher}{Society of Exploration Geophysicists}. pp.
  \bibinfo{pages}{81--82}.
\bibitem[{Constable et~al.(1987)Constable, Parker and
  Constable}]{constable1987occam}
\bibinfo{author}{Constable, S.C.}, \bibinfo{author}{Parker, R.L.},
  \bibinfo{author}{Constable, C.G.}, \bibinfo{year}{1987}.
\newblock \bibinfo{title}{{Occam’s inversion: A practical algorithm for
  generating smooth models from electromagnetic sounding data}}.
\newblock \bibinfo{journal}{Geophysics} \bibinfo{volume}{52},
  \bibinfo{pages}{289--300}.
\bibitem[{Eidesmo et~al.(2002)Eidesmo, Ellingsrud, MacGregor, Constable, Sinha,
  Johansen, Kong and Westerdahl}]{eidesmo2002sbl}
\bibinfo{author}{Eidesmo, T.}, \bibinfo{author}{Ellingsrud, S.},
  \bibinfo{author}{MacGregor, L.}, \bibinfo{author}{Constable, S.},
  \bibinfo{author}{Sinha, M.}, \bibinfo{author}{Johansen, S.},
  \bibinfo{author}{Kong, F.}, \bibinfo{author}{Westerdahl, H.},
  \bibinfo{year}{2002}.
\newblock \bibinfo{title}{{Sea bed logging (SBL), a new method for remote and
  direct identification of hydrocarbon filled layers in deepwater areas}}.
\newblock \bibinfo{journal}{First break} \bibinfo{volume}{20}.
\bibitem[{Ellingsrud et~al.(2002)Ellingsrud, Eidesmo, Johansen, Sinha,
  MacGregor and Constable}]{ellingsrud2002remote}
\bibinfo{author}{Ellingsrud, S.}, \bibinfo{author}{Eidesmo, T.},
  \bibinfo{author}{Johansen, S.}, \bibinfo{author}{Sinha, M.},
  \bibinfo{author}{MacGregor, L.}, \bibinfo{author}{Constable, S.},
  \bibinfo{year}{2002}.
\newblock \bibinfo{title}{{Remote sensing of hydrocarbon layers by seabed
  logging (SBL): Results from a cruise offshore Angola}}.
\newblock \bibinfo{journal}{The Leading Edge} \bibinfo{volume}{21},
  \bibinfo{pages}{972--982}.
\bibitem[{Grayver et~al.(2013)Grayver, Streich and Ritter}]{grayver2013gji}
\bibinfo{author}{Grayver, A.V.}, \bibinfo{author}{Streich, R.},
  \bibinfo{author}{Ritter, O.}, \bibinfo{year}{2013}.
\newblock \bibinfo{title}{{Three-dimensional parallel distributed inversion of
  CSEM data using a direct forward solver}}.
\newblock \bibinfo{journal}{Geophysical Journal International}
  \bibinfo{volume}{193}, \bibinfo{pages}{1432--1446}.
\bibitem[{Grayver et~al.(2014)Grayver, Streich and Ritter}]{grayver2014geo}
\bibinfo{author}{Grayver, A.V.}, \bibinfo{author}{Streich, R.},
  \bibinfo{author}{Ritter, O.}, \bibinfo{year}{2014}.
\newblock \bibinfo{title}{{3D inversion and resolution analysis of land-based
  CSEM data from the Ketzin CO2 storage formation}}.
\newblock \bibinfo{journal}{Geophysics} \bibinfo{volume}{79},
  \bibinfo{pages}{E101--E114}.
\bibitem[{Gribenko and Zhdanov(2007)}]{gribenko2007rigorous}
\bibinfo{author}{Gribenko, A.}, \bibinfo{author}{Zhdanov, M.},
  \bibinfo{year}{2007}.
\newblock \bibinfo{title}{{Rigorous 3D inversion of marine CSEM data based on
  the integral equation method}}.
\newblock \bibinfo{journal}{Geophysics} \bibinfo{volume}{72},
  \bibinfo{pages}{WA73--WA84}.
\bibitem[{Key(2016)}]{key2016mare2dem}
\bibinfo{author}{Key, K.}, \bibinfo{year}{2016}.
\newblock \bibinfo{title}{{MARE2DEM: a 2-D inversion code for controlled-source
  electromagnetic and magnetotelluric data}}.
\newblock \bibinfo{journal}{Geophysical Journal International}
  \bibinfo{volume}{207}, \bibinfo{pages}{571--588}.
\bibitem[{Komatitsch and Martin(2007)}]{Komatitsch_2007_GEO}
\bibinfo{author}{Komatitsch, D.}, \bibinfo{author}{Martin, R.},
  \bibinfo{year}{2007}.
\newblock \bibinfo{title}{{An unsplit convolutional perfectly matched layer
  improved at grazing incidence for the seismic wave equation}}.
\newblock \bibinfo{journal}{Geophysics} \bibinfo{volume}{72},
  \bibinfo{pages}{SM155--SM167}.
\bibitem[{Lee et~al.(1989)Lee, Liu and Morrison}]{lee1989new}
\bibinfo{author}{Lee, K.H.}, \bibinfo{author}{Liu, G.},
  \bibinfo{author}{Morrison, H.}, \bibinfo{year}{1989}.
\newblock \bibinfo{title}{A new approach to modeling the electromagnetic
  response of conductive media}.
\newblock \bibinfo{journal}{Geophysics} \bibinfo{volume}{54},
  \bibinfo{pages}{1180--1192}.
\bibitem[{Lelièvre and Farquharson(2013)}]{Lelievre2013gradient}
\bibinfo{author}{Lelièvre, P.G.}, \bibinfo{author}{Farquharson, C.G.},
  \bibinfo{year}{2013}.
\newblock \bibinfo{title}{{Gradient and smoothness regularization operators for
  geophysical inversion on unstructured meshes}}.
\newblock \bibinfo{journal}{Geophysical Journal International}
  \bibinfo{volume}{195}, \bibinfo{pages}{330--341}.
\newblock \URLprefix \url{https://doi.org/10.1093/gji/ggt255},
  \DOIprefix\doi{10.1093/gji/ggt255}.
\bibitem[{Li and Key(2007)}]{li20072d}
\bibinfo{author}{Li, Y.}, \bibinfo{author}{Key, K.}, \bibinfo{year}{2007}.
\newblock \bibinfo{title}{{2D} marine controlled-source electromagnetic
  modeling: Part 1—an adaptive finite-element algorithm}.
\newblock \bibinfo{journal}{Geophysics} \bibinfo{volume}{72},
  \bibinfo{pages}{WA51--WA62}.
\bibitem[{Maa{\o}(2007)}]{Maao_2007_FFT}
\bibinfo{author}{Maa{\o}, F.}, \bibinfo{year}{2007}.
\newblock \bibinfo{title}{Fast finite-difference time-domain modeling for
  marine subsurface electromagnetic problems}.
\newblock \bibinfo{journal}{Geophysics} \bibinfo{volume}{72},
  \bibinfo{pages}{A19--A23}.
\bibitem[{MacGregor et~al.(2007)MacGregor, Barker, Overton, Moody and
  Bodecott}]{macgregor2007derisking}
\bibinfo{author}{MacGregor, L.}, \bibinfo{author}{Barker, N.},
  \bibinfo{author}{Overton, A.}, \bibinfo{author}{Moody, S.},
  \bibinfo{author}{Bodecott, D.}, \bibinfo{year}{2007}.
\newblock \bibinfo{title}{{Derisking exploration prospects using integrated
  seismic and electromagnetic data-A Falkland Islands case study}}.
\newblock \bibinfo{journal}{The Leading Edge} \bibinfo{volume}{26},
  \bibinfo{pages}{356--359}.
\bibitem[{MacGregor and Tomlinson(2014)}]{macgregor2014mcsem}
\bibinfo{author}{MacGregor, L.}, \bibinfo{author}{Tomlinson, J.},
  \bibinfo{year}{2014}.
\newblock \bibinfo{title}{{Marine controlled-source electromagnetic methods in
  the hydrocarbon industry: A tutorial on method and practice}}.
\newblock \bibinfo{journal}{Interpretation} \bibinfo{volume}{2},
  \bibinfo{pages}{SH13--SH32}.
\bibitem[{Mittet(2010)}]{Mittet_2010_HFD}
\bibinfo{author}{Mittet, R.}, \bibinfo{year}{2010}.
\newblock \bibinfo{title}{High-order finite-difference simulations of marine
  {CSEM} surveys using a correspondence principle for wave and diffusion
  fields}.
\newblock \bibinfo{journal}{Geophysics} \bibinfo{volume}{75},
  \bibinfo{pages}{F33--F50}.
\bibitem[{Mittet and Morten(2012)}]{mittet2012detection}
\bibinfo{author}{Mittet, R.}, \bibinfo{author}{Morten, J.P.},
  \bibinfo{year}{2012}.
\newblock \bibinfo{title}{Detection and imaging sensitivity of the marine
  {CSEM} method}.
\newblock \bibinfo{journal}{Geophysics} \bibinfo{volume}{77},
  \bibinfo{pages}{E411--E425}.
\bibitem[{Morten et~al.(2009)Morten, Bj{\o}rke and
  St{\o}ren}]{morten2009uncertainy}
\bibinfo{author}{Morten, J.P.}, \bibinfo{author}{Bj{\o}rke, A.K.},
  \bibinfo{author}{St{\o}ren, T.}, \bibinfo{year}{2009}.
\newblock \bibinfo{title}{{CSEM data uncertainty analysis for 3D inversion}},
  in: \bibinfo{booktitle}{SEG Technical Program Expanded Abstracts 2009},
  \bibinfo{publisher}{Society of Exploration Geophysicists}. pp.
  \bibinfo{pages}{724--728}.
\bibitem[{Mulder(2006)}]{mulder2006multigrid}
\bibinfo{author}{Mulder, W.}, \bibinfo{year}{2006}.
\newblock \bibinfo{title}{A multigrid solver for 3{D} electromagnetic
  diffusion}.
\newblock \bibinfo{journal}{Geophysical prospecting} \bibinfo{volume}{54},
  \bibinfo{pages}{633--649}.
\bibitem[{Newman and Alumbaugh(1995)}]{newman1995frequency}
\bibinfo{author}{Newman, G.A.}, \bibinfo{author}{Alumbaugh, D.L.},
  \bibinfo{year}{1995}.
\newblock \bibinfo{title}{Frequency-domain modelling of airborne
  electromagnetic responses using staggered finite differences}.
\newblock \bibinfo{journal}{Geophysical Prospecting} \bibinfo{volume}{43},
  \bibinfo{pages}{1021--1042}.
\bibitem[{Nocedal and Wright(2006)}]{Nocedal_2006_NOO}
\bibinfo{author}{Nocedal, J.}, \bibinfo{author}{Wright, S.J.},
  \bibinfo{year}{2006}.
\newblock \bibinfo{title}{Numerical Optimization}.
\newblock \bibinfo{edition}{2nd} ed., \bibinfo{publisher}{Springer}.
\bibitem[{Oristaglio and Hohmann(1984)}]{oristaglio1984diffusion}
\bibinfo{author}{Oristaglio, M.L.}, \bibinfo{author}{Hohmann, G.W.},
  \bibinfo{year}{1984}.
\newblock \bibinfo{title}{Diffusion of electromagnetic fields into a
  two-dimensional earth: A finite-difference approach}.
\newblock \bibinfo{journal}{Geophysics} \bibinfo{volume}{49},
  \bibinfo{pages}{870--894}.
\bibitem[{Paige and Saunders(1982)}]{Paige_1982_ALS}
\bibinfo{author}{Paige, C.C.}, \bibinfo{author}{Saunders, M.A.},
  \bibinfo{year}{1982}.
\newblock \bibinfo{title}{{ALGORITHM 583} {LSQR} : Sparse linear equations and
  least squares problems}.
\newblock \bibinfo{journal}{{ACM} Transactions on Mathematical Software}
  \bibinfo{volume}{8}, \bibinfo{pages}{195--209}.
\bibitem[{Plessix and Mulder(2008)}]{Plessix_2008_RIC}
\bibinfo{author}{Plessix, R.E.}, \bibinfo{author}{Mulder, W.A.},
  \bibinfo{year}{2008}.
\newblock \bibinfo{title}{Resistivity imaging with controlled-source
  electromagnetic data: depth and data weighting}.
\newblock \bibinfo{journal}{Inverse Problems} \bibinfo{volume}{24},
  \bibinfo{pages}{034012}.
\bibitem[{Puzyrev et~al.(2013)Puzyrev, Koldan, de~la Puente, Houzeaux,
  V{\'a}zquez and Cela}]{puzyrev2013parallel}
\bibinfo{author}{Puzyrev, V.}, \bibinfo{author}{Koldan, J.},
  \bibinfo{author}{de~la Puente, J.}, \bibinfo{author}{Houzeaux, G.},
  \bibinfo{author}{V{\'a}zquez, M.}, \bibinfo{author}{Cela, J.M.},
  \bibinfo{year}{2013}.
\newblock \bibinfo{title}{A parallel finite-element method for
  three-dimensional controlled-source electromagnetic forward modelling}.
\newblock \bibinfo{journal}{Geophysical Journal International}
  \bibinfo{volume}{193}, \bibinfo{pages}{678--693}.
\bibitem[{Rochlitz et~al.(2019)Rochlitz, Skibbe and
  G{\"u}nther}]{rochlitz2019custem}
\bibinfo{author}{Rochlitz, R.}, \bibinfo{author}{Skibbe, N.},
  \bibinfo{author}{G{\"u}nther, T.}, \bibinfo{year}{2019}.
\newblock \bibinfo{title}{{custEM: Customizable finite-element simulation of
  complex controlled-source electromagnetic data}}.
\newblock \bibinfo{journal}{Geophysics} \bibinfo{volume}{84},
  \bibinfo{pages}{F17--F33}.
\bibitem[{Saad(2003)}]{Saad_2003_IMS}
\bibinfo{author}{Saad, Y.}, \bibinfo{year}{2003}.
\newblock \bibinfo{title}{Iterative {M}ethods for {S}parse {L}inear {S}ystems}.
\newblock \bibinfo{publisher}{SIAM}, \bibinfo{address}{Philadelphia}.
\bibitem[{Schwarzbach and Haber(2013)}]{schwarzbach2013gji}
\bibinfo{author}{Schwarzbach, C.}, \bibinfo{author}{Haber, E.},
  \bibinfo{year}{2013}.
\newblock \bibinfo{title}{Finite element based inversion for time-harmonic
  electromagnetic problems}.
\newblock \bibinfo{journal}{Geophysical Journal International}
  \bibinfo{volume}{193}, \bibinfo{pages}{615--634}.
\bibitem[{Shantsev et~al.(2020)Shantsev, Nerland and
  Gelius}]{shantsev2002timelapse}
\bibinfo{author}{Shantsev, D.V.}, \bibinfo{author}{Nerland, E.A.},
  \bibinfo{author}{Gelius, L.J.}, \bibinfo{year}{2020}.
\newblock \bibinfo{title}{{Time-lapse CSEM: how important is survey
  repeatability?}}
\newblock \bibinfo{journal}{Geophysical Journal International}
  \bibinfo{volume}{223}, \bibinfo{pages}{2133--2147}.
\bibitem[{da~Silva et~al.(2012)da~Silva, Morgan, MacGregor and
  Warner}]{da2012finite}
\bibinfo{author}{da~Silva, N.V.}, \bibinfo{author}{Morgan, J.V.},
  \bibinfo{author}{MacGregor, L.}, \bibinfo{author}{Warner, M.},
  \bibinfo{year}{2012}.
\newblock \bibinfo{title}{{A finite element multifrontal method for {3D} CSEM
  modeling in the frequency domain}}.
\newblock \bibinfo{journal}{Geophysics} \bibinfo{volume}{77},
  \bibinfo{pages}{E101--E115}.
\bibitem[{Sirgue et~al.(2008)Sirgue, Etgen and Albertin}]{Sirgue_2008_FDW}
\bibinfo{author}{Sirgue, L.}, \bibinfo{author}{Etgen, J.T.},
  \bibinfo{author}{Albertin, U.}, \bibinfo{year}{2008}.
\newblock \bibinfo{title}{3{D} {F}requency {D}omain {W}aveform {I}nversion
  using {T}ime {D}omain {F}inite {D}ifference {M}ethods}, in:
  \bibinfo{booktitle}{Proceedings 70th {EAGE}, Conference and Exhibition, Roma,
  Italy}, p. \bibinfo{pages}{F022}.
\bibitem[{Sirgue et~al.(2010)Sirgue, Etgen, Albertin and
  Brandsberg-Dahl}]{sirgue2010system}
\bibinfo{author}{Sirgue, L.}, \bibinfo{author}{Etgen, J.T.},
  \bibinfo{author}{Albertin, U.}, \bibinfo{author}{Brandsberg-Dahl, S.},
  \bibinfo{year}{2010}.
\newblock \bibinfo{title}{System and method for 3{D} frequency domain waveform
  inversion based on 3{D} time-domain forward modeling}.
\newblock \bibinfo{note}{US Patent 7,725,266}.
\bibitem[{Smith(1996a)}]{smith1996conservative1}
\bibinfo{author}{Smith, J.T.}, \bibinfo{year}{1996}a.
\newblock \bibinfo{title}{Conservative modeling of 3-{D} electromagnetic
  fields, {P}art i: Properties and error analysis}.
\newblock \bibinfo{journal}{Geophysics} \bibinfo{volume}{61},
  \bibinfo{pages}{1308--1318}.
\bibitem[{Smith(1996b)}]{smith1996conservative2}
\bibinfo{author}{Smith, J.T.}, \bibinfo{year}{1996}b.
\newblock \bibinfo{title}{{Conservative modeling of 3-D electromagnetic fields,
  Part II: Biconjugate gradient solution and an accelerator}}.
\newblock \bibinfo{journal}{Geophysics} \bibinfo{volume}{61},
  \bibinfo{pages}{1319--1324}.
\bibitem[{St{\o}ren et~al.(2008)St{\o}ren, Zach and
  Maa{\o}}]{Storen_2008_Gradient}
\bibinfo{author}{St{\o}ren, T.}, \bibinfo{author}{Zach, J.J.},
  \bibinfo{author}{Maa{\o}, F.A.}, \bibinfo{year}{2008}.
\newblock \bibinfo{title}{Gradient calculations for 3{D} inversion of {CSEM}
  data using a fast finite-difference time-domain modelling code}, in:
  \bibinfo{booktitle}{70th EAGE Conference and Exhibition incorporating SPE
  EUROPEC 2008}.
\bibitem[{Streich(2009)}]{streich20093d}
\bibinfo{author}{Streich, R.}, \bibinfo{year}{2009}.
\newblock \bibinfo{title}{3{D} finite-difference frequency-domain modeling of
  controlled-source electromagnetic data: {D}irect solution and optimization
  for high accuracy}.
\newblock \bibinfo{journal}{Geophysics} \bibinfo{volume}{74},
  \bibinfo{pages}{F95--F105}.
\bibitem[{Symes(2007)}]{Symes_2007_RTM}
\bibinfo{author}{Symes, W.W.}, \bibinfo{year}{2007}.
\newblock \bibinfo{title}{Reverse time migration with optimal checkpointing}.
\newblock \bibinfo{journal}{Geophysics} \bibinfo{volume}{72},
  \bibinfo{pages}{SM213--SM221}.
\newblock \URLprefix \url{http://link.aip.org/link/?GPY/72/SM213/1},
  \DOIprefix\doi{10.1190/1.2742686}.
\bibitem[{Taflove and Hagness(2005)}]{Taflove_2005_CEF}
\bibinfo{author}{Taflove, A.}, \bibinfo{author}{Hagness, S.C.},
  \bibinfo{year}{2005}.
\newblock \bibinfo{title}{Computational Electrodynamics: The Finite-Difference
  Time-Domain Method}.
\newblock \bibinfo{edition}{3rd} ed., \bibinfo{publisher}{Artech House}.
\bibitem[{Wang and Hohmann(1993)}]{wang1993finite}
\bibinfo{author}{Wang, T.}, \bibinfo{author}{Hohmann, G.W.},
  \bibinfo{year}{1993}.
\newblock \bibinfo{title}{A finite-difference, time-domain solution for
  three-dimensional electromagnetic modeling}.
\newblock \bibinfo{journal}{Geophysics} \bibinfo{volume}{58},
  \bibinfo{pages}{797--809}.
\bibitem[{Ward and Hohmann(1988)}]{ward1988electromagnetic}
\bibinfo{author}{Ward, S.H.}, \bibinfo{author}{Hohmann, G.W.},
  \bibinfo{year}{1988}.
\newblock \bibinfo{title}{Electromagnetic theory for geophysical applications},
  in: \bibinfo{booktitle}{{Electromagnetic Methods in Applied Geophysics:
  Volume 1, Theory}}. \bibinfo{publisher}{Society of Exploration
  Geophysicists}, pp. \bibinfo{pages}{130--311}.
\bibitem[{Yang(2023)}]{Yang_2023_libEMM}
\bibinfo{author}{Yang, P.}, \bibinfo{year}{2023}.
\newblock \bibinfo{title}{{libEMM: A fictious wave domain 3D CSEM modelling
  library bridging sequential and parallel GPU implementation}}.
\newblock \bibinfo{journal}{Computer Physics Communications}
  \bibinfo{volume}{288}, \bibinfo{pages}{108745}.
\newblock \DOIprefix\doi{https://doi.org/10.1016/j.cpc.2023.108745}.
\bibitem[{Yang et~al.(2016a)Yang, Brossier, M\'etivier and
  Virieux}]{Yang_2016_CAR}
\bibinfo{author}{Yang, P.}, \bibinfo{author}{Brossier, R.},
  \bibinfo{author}{M\'etivier, L.}, \bibinfo{author}{Virieux, J.},
  \bibinfo{year}{2016}a.
\newblock \bibinfo{title}{Wavefield reconstruction in attenuating media: A
  checkpointing-assisted reverse-forward simulation method}.
\newblock \bibinfo{journal}{Geophysics} \bibinfo{volume}{81},
  \bibinfo{pages}{R349--R362}.
\newblock \DOIprefix\doi{10.1190/geo2016-0082.1}.
\bibitem[{Yang et~al.(2016b)Yang, Brossier and Virieux}]{Yang_2016_WRB}
\bibinfo{author}{Yang, P.}, \bibinfo{author}{Brossier, R.},
  \bibinfo{author}{Virieux, J.}, \bibinfo{year}{2016}b.
\newblock \bibinfo{title}{Wavefield reconstruction from significantly decimated
  boundaries}.
\newblock \bibinfo{journal}{Geophysics} \bibinfo{volume}{80},
  \bibinfo{pages}{T197--T209}.
\newblock \DOIprefix\doi{10.1190/GEO2015-0711.1}.
\bibitem[{Yang and Mittet(2023)}]{Yang_2023_HFDNU}
\bibinfo{author}{Yang, P.}, \bibinfo{author}{Mittet, R.}, \bibinfo{year}{2023}.
\newblock \bibinfo{title}{Controlled-source electromagnetics modelling using
  high order finite-difference time-domain method on a nonuniform grid}.
\newblock \bibinfo{journal}{Geophysics} \bibinfo{volume}{88},
  \bibinfo{pages}{E53--E67}.
\newblock \DOIprefix\doi{10.1190/geo2022-0134.1}.
\bibitem[{Zaslavsky et~al.(2013)Zaslavsky, Druskin, Abubakar, Habashy and
  Simoncini}]{zaslavsky2013large}
\bibinfo{author}{Zaslavsky, M.}, \bibinfo{author}{Druskin, V.},
  \bibinfo{author}{Abubakar, A.}, \bibinfo{author}{Habashy, T.},
  \bibinfo{author}{Simoncini, V.}, \bibinfo{year}{2013}.
\newblock \bibinfo{title}{Large-scale gauss-newton inversion of transient csem
  data using the model order reduction framework}.
\newblock \bibinfo{journal}{Geophysics} \bibinfo{volume}{78},
  \bibinfo{pages}{E161--E171}.
\bibitem[{Zhdanov and Keller(1994)}]{zhdanov1994geoelectrical}
\bibinfo{author}{Zhdanov, M.S.}, \bibinfo{author}{Keller, G.V.},
  \bibinfo{year}{1994}.
\newblock \bibinfo{title}{The geoelectrical methods in geophysical
  exploration}. volume~\bibinfo{volume}{31}.

\end{thebibliography}

\end{document}